\documentclass[a4paper]{article}

\usepackage{amsmath,amssymb,amsfonts}%
\usepackage{amsthm,mathtools}%

\usepackage{authblk}%
\usepackage[numbers]{natbib}%
\usepackage[colorlinks=false,pdfborder={0 0 0}]{hyperref}%
\usepackage{enumitem}%

\usepackage{threeparttable}
\usepackage{booktabs}
\usepackage{multirow}
\usepackage{adjustbox}

\usepackage{orcidlink}
\usepackage{doi}

\newcommand{\N}{\mathbb{N}}
\newcommand{\Z}{\mathbb{Z}}
\newcommand{\R}{\mathbb{R}}
\newcommand{\Q}{\mathbb{Q}}

\newcommand{\poly}{\text{P}}
\newcommand{\npoly}{\text{NP}}
\newcommand{\conpoly}{\text{coNP}}

\newcommand{\gdw}{\Leftrightarrow}
\newcommand{\st}{\text{s.t.}}

\DeclareMathOperator{\rank}{rank}
\DeclareMathOperator{\conv}{conv}

\DeclareMathOperator{\face}{\mathcal F}

\newcommand{\e}{\mathbf{e}}

\newtheorem{theorem}{Theorem}%
\newtheorem{lemma}[theorem]{Lemma}
\newtheorem{corollary}{Corollary}

\newtheorem{remark}{Remark}%

\newtheorem{definition}{Definition}%

\begin{document}

\title{On the feasibility of generalized\\ inverse linear programs}

\author{Christoph Buchheim~\orcidlink{0000-0001-9974-404X} and Lowig T.\ Duer~\orcidlink{0009-0002-1990-4029}}

\affil{Fakultät für Mathematik, TU Dortmund, Germany}

\maketitle

\begin{abstract}
    We investigate the feasibility problem for generalized inverse linear programs. Given an LP with affinely parametrized objective function and right-hand side as well as a target set~$Y$, the goal is to decide whether the parameters can be chosen such that there exists an optimal solution that belongs to~$Y$ (optimistic scenario) or such that all optimal solutions belong to~$Y$ (pessimistic scenario).
    We study the complexity of this decision problem and show how it depends on the structure of the set~$Y$, the form of the LP, the adjustable parameters, and the underlying scenario.
    For a target singleton~$Y = \{\bar y\}$, we show that the problem is tractable if the given LP is in standard form, but $\npoly$-hard if the LP is given in natural form.
    If instead we are given a target basis~$\bar{\mathcal B}$, the problem in standard form becomes $\npoly$-complete in the optimistic case, while remaining tractable in the pessimistic case.
    For partially fixed target solutions, the problem gets almost immediately $\npoly$-hard, but for particular cases we prove tractability if the number of free target variables is fixed.
    Moreover, we give a rigorous proof of membership in $\npoly$ for any polyhedral target set, and discuss how this property can be extended to more general target sets using an oracle-based approach.\\[1ex]
\end{abstract}

\section{Introduction}\label{sec:intro}

Inverse optimization is a fundamental problem in optimization that has been studied extensively in the past.
Contrary to a standard forward optimization problem, where the goal is to find a solution that optimizes a given objective function while satisfying some given constraints, an inverse optimization problem arises when a solution to a forward problem is already known but instead some parameters of the problem need to be adjusted such that the given solution becomes optimal.
Inverse linear optimization describes such inverse problems when the initial forward problem is a linear program (LP).
Usually, the parameters to be adjusted are the objective function coefficients, and the optimization goal of the inverse problem is to minimize the perturbation of the initially given objective function with respect to some norm.
In partial inverse problems, the desired optimal solution is not completely known, but only partially, e.g., only some variables are fixed.
Then, the goal is to adjust the parameters of the forward problem such that there exists an optimal extension of the partial solution.
Therefore, partial inverse problems form a generalization of inverse problems.

Typical applications of inverse optimization are found in estimation problems and incentive design.
Estimation problems arise when optimization problems are used to model any kind of system, but some parameters of the model are difficult to determine.
With inverse optimization, prior observation can be used to determine plausible parameters for a fitting model, e.g., by considering past route data of experienced drivers to improve the coefficients of a capacitated vehicle routing problem \cite{Chen2021inverseoptapproach} or by using openly accessible energy market 
data to reconstruct non-public offer prices~
\cite{Liang2023data-driven}.
Incentive design, on the other hand, describes a field where the goal is to indirectly incentivize someone else's decision making process towards a desired behavior by setting some influencing parameters.
Many such applications arise when an authority wants to elicit a certain behavior of other players like citizens or companies by enacting laws and regulations and anticipating their consequences on the other player's behavior.
E.g., in toll pricing policies, inverse linear programs can model the problem of an authority to set toll prices on existing road networks such that users choose ``system-optimizing'' routes 
with minimum average time \cite{Dial1999minimal-revenue}, minimum environmental impact \cite{Li2019environmentalToll}, or risk-minimization for the transport of hazardous materials \cite{Marcotte2009tollpolicies}.
Similarly, in climate politics, governmental incentives to improve private investments in renewable energy sources can be modeled with inverse optimization, as Zhou et al.~\cite{Zhou2011renewableenergy} do by decomposing a bilevel problem into a partial inverse MILP.
The field of incentive design illustrates in particular the relation of inverse optimization and bilevel optimization, which will be further elaborated in the following paragraph.

In many applications, it is reasonable to assume that there exists more than one target solution. For instance, in many estimation problems the outcome of the parametrized optimization can only be observed partially or only up to a certain precision. In incentive design, it often suffices that the induced decisions lie within a certain range. 
This motivates our investigation of generalized inverse problems.

\vspace*{-2ex}
\paragraph{Related literature.}
Inverse optimization seems to have been first studied by geophysicists in the 1970s; see~\cite{tarantola1987inverseproblemtheory}.
In the 1990s, some research has been conducted on the inverse variants of several combinatorial problems \cite{Burton1992inverseshortestpath,Zhang1996inverseminspan}, before in 1996 the inverse linear problem was first formulated and proven to be tractable under the~$\ell_1$- and~$\ell_\infty$-norm by reducing it to an LP \cite{Zhang1996invlinprog}.
Later, this complexity result was famously extended to inverse problems of arbitrary tractable forward problems with linear objective functions \cite{ahuja_orlin1995inverse_optimization} using fundamental complexity results based on the ellipsoid method by Grötschel et al.~\cite{grotschel1993geometricalgorithms}.
More recently, it has been shown for inverse MILPs 
that generally the primal bound decision variant is $\conpoly$-complete, while the dual bound decision variant is $\npoly$-complete, conversely to the respective forward problems~\cite{Bulut2021InvMILP}.
There exist multiple surveys on inverse optimization, of which we want to mention two; a survey from Heuberger~\cite{heuberger2004survey} with a focus on combinatorial problems and a newer survey by Chan et al.~\cite{Chan2023Survey} with a more general summary of solution techniques, theory, and applications.

In recent years, partial inverse optimization has become a subject of intensive research as well.
Most publications in this field consider partial inverse versions of specific combinatorial forward problems where only part of the variables are fixed for the optimal solution.
In contrast to standard inverse optimization, the complexity of partial inverse problems can increase relative to the forward problem.
The assignment problem is the first such tractable problem for which the partial inverse version has been shown to be $\npoly$-hard if the change of coefficients is bounded \cite{Yang2001PartInvAssgn}, whereas for unbounded coefficients the complexity differs depending on whether the $\ell_1$- or $\ell_\infty$-norm is considered~\cite{Yang2007ParInvAssgn,lai2003complexity}.
In several publications, a distinction is made between variants where the values of the objective function may only be decreased or increased, respectively, 
as 
the first variant 
tends 
to be easier than the latter.
Very recently, making this distinction, Ley and Merkert~\cite{Ley2025solution} have studied several partial inverse combinatorial problems, summarized prior results, and provided a tabular overview of results to date.
They also proved that for a partial inverse shortest path problem on directed graphs, it is even $\npoly$-hard to decide if an instance is feasible at all.
For partial inverse LPs, there exists a report by Gentry~\cite{gentry2001partial} discussing MILP reformulations~\mbox{(see \cite{Yang2007ParInvAssgn})}, but it is not publicly available.
Generally, however, there seems to be less research on partial inverse linear problems and more general classes of inverse linear problems so far.

Kurtz et al.~\cite{Kurtz2026CounterfactualLP} study a problem that is closely related to inverse optimization called \emph{counterfactual explanations}. They ask
which values certain input parameters of an optimization problem should have 
in order for the optimal solution to be part of some favored solution space, which corresponds to our target set.
Their \emph{weak} and \emph{strong} problems correspond to our optimistic and pessimistic scenario.
Among other things, they show that even deciding if a counterfactual explanation exists can be $\npoly$-hard when the favored solution space is a polyhedron.
In contrast to our model, they also consider the constraints to be adjustable input parameters (in addition to objective function and right-hand side), thus our lower bounds on the complexity will be stronger results.
Furthermore, Kurtz et al.~\cite{Kurtz2025CounterfactualILP} perform a very similar analysis on counterfactual explanations for (possibly) $\npoly$-hard integer optimization problems.
They prove completeness for $\Sigma_2^\poly$ for their resulting weak and strong problems and even specify that the hardness is preserved when certain input parameters remain fixed.
Algorithms to compute counterfactual explanations have also been developed, e.g., for LPs based on a single-level reformulation of the problem~\cite{lefebvre2024counterfactual} and for MILPs~\cite{Korikov2023Counterfactual}.

Another closely related field to (partial) inverse optimization is bilevel optimization, a class of optimization problems of game theoretical nature, where two players make consecutive decisions that depend on each other. For further details on bilevel optimization, we refer to~\cite{Colson2007Bileveloverview,Dempe2015BilevelTheorAlgoAppl}.
Many inverse problems can be formulated as bilevel programs, where the lower level consist of the forward problem and the upper level contains the minimization of the parameter perturbation and the inverse optimality condition.
Vice versa, asking whether there exists an upper-level solution that enforces a specific 
lower-level response corresponds to the question of feasibility of an inverse program.
Investigating the complexity of bilevel spanning tree problems, Buchheim et al.~\cite{BuchheimHenke2022BilevelMST} study this question and show that for bilevel minimum spanning tree it is $\npoly$-hard to answer, even if only one variable is fixed in the partial solution.
Bilevel optimization is well known to be~$\npoly$-hard~\cite{jeroslow1985bilevel}, 
even with a single upper-level variable~\cite{Sugishita2025Bilevelwithsinglevariable}.

\vspace*{-2ex}
\paragraph{Our contribution.}
We investigate a generalized version of (partial) inverse linear optimization, in particular the complexity of deciding feasibility.
Instead of requiring some coefficients of the optimal solution to attain given values, we require that the optimal solution belongs to some given \emph{target set}.
So starting with a parametrized LP as a forward problem and a target set, the goal is to decide whether one element of the target set can be made optimal by adjusting the parameters of the program.
We will study the complexity for different such target sets, mainly considering special classes of polyhedra.

However, the addressed inverse linear problems will also be more general than many common definitions of (partial) inverse programs in some other aspects:
We assume that not only the objective function but also the right-hand sides of the linear constraints can be adjusted.
Moreover, both adjustable parameters might not be directly selectable but are determined by other parameters through affine linear functions.
As a result, both parameters can be coupled, so that possible dependencies can be modeled.
We will limit our study to the problem of deciding feasibility, which already turns out to be a complex task with surprising results in our general setting.
Compared to existing complexity results (especially from~\cite{Kurtz2026CounterfactualLP}) we perform a much more detailed analysis on which aspects of the inverse problem actually determine the complexity and how they need to be restricted to reduce the complexity.
Moreover, beyond hardness, we also show completeness results. 

In fact, we will first focus on the special case where a single complete target solution is given (corresponding to a standard inverse linear problem) and show that, in this case, deciding feasibility is tractable if the forward LP is given in standard form, but $\npoly$-complete if it is given in natural form.
If instead of a complete solution we ask for a basis to become optimal, the problem is $\npoly$-complete in the optimistic scenario and tractable in the pessimistic scenario.
If the target solution is only given partially and not all variables are fixed, the problem in standard form becomes $\npoly$-complete as well, even if only one variable is fixed.
However, we can show fixed-parameter tractability in the number of free variables of the partial solution.
We also prove for several cases that the problem can be easier if only the objective function is parametrized.
Moreover, we show that in all considered cases the generalized inverse feasibility problem belongs to $\npoly$, i.e., feasibility can be verified efficiently with an appropriate certificate.
For the optimistic scenario, we show that (to a certain extent) this can even be generalized to instances where the target set is only provided via some type of weak membership oracle.

\vspace*{-2ex}
\paragraph{Outline.}
The remainder of this paper is organized as follows.
In Section~\ref{sec:prel}, we present a precise formulation of our \emph{generalized inverse linear feasibility problem}, specify some notation, and recall necessary preliminaries on LP optimality conditions.
Next, we study the complexity of the inverse linear feasibility problem for a single target solution~(Section~\ref{sec:ciflp}) and for a given target basis~(Section~\ref{sec:targetbasis}).
In Section~\ref{sec:piflp}, we consider polyhedral target sets, including the case of a partially given target solution.
In Section~\ref{sec:oracle-target}, we study the problem variant where the target set is only provided by a weak membership oracle. Section~\ref{sec:conclusion} concludes.

\section{Preliminaries}\label{sec:prel}

In our general setting, we consider a parametric LP in natural form
\begin{equation}
    \left. \begin{aligned}
        \min_{y}&       \quad   (Cx+c)^\top y \\
        \st&            \quad   Ay \leq Bx+b,\ y\in\R^n
    \end{aligned}\qquad \right\} \tag{LP$_x$} \label{p:lpx}
\end{equation}
where, in addition to the constant parameters $A\in\R^{m\times n},\ b\in\R^m$, and~$c\in\R^n$, the objective function and right hand side depend linearly on a parameter~\mbox{$x\in\R^k$} via matrices $B\in\R^{m\times k}$ and $C\in\R^{n\times k}$.
In some (clearly indicated) cases we will consider a similar parametric LP in standard form
\begin{equation}
    \left. \begin{aligned}
        \min_{y}&       \quad   (Cx+c)^\top y \\
        \st&            \quad   Ay = Bx+b,\ y\in\R^n_+\;,
    \end{aligned}\qquad \right\} \tag{LP$'_x$} \label{p:lpx-stand}
\end{equation}
where $\R_+^n\coloneqq \{y\in\R^n \mid y\geq 0\}$.
In both cases, we assume that the parameter~$x$ is restricted to a polyhedron~$X\subseteq\R^k$ called \emph{parameter set} and that a \emph{target set}~$Y\subseteq\R^n$ in the variable space of the LP is given.

We are now interested in the following question: 
Does there exist an~$x\in X$ such that the optimal solution to (\hyperref[p:lpx]{LP$_x$}) belongs to $Y$?
We will call this decision problem the \emph{(Generalized) Inverse Linear Feasibility Problem} \hypertarget{p:ilfp}{(InvLFP)}.
For a formal definition of this problem, we need to take into account that optimal solutions are not necessarily unique, so that we have to distinguish between the \emph{optimistic} and \emph{pessimistic} scenario.
This distinction is less common in inverse optimization, where usually the optimistic scenario is considered, but rather motivated by bilevel optimization, where we assume that the parameter $x$ and the optimization variable $y$ are selected by different players.
The two scenarios can be interpreted as the behavior of a follower 
solving the forward LP seen from the perspective of a leader
choosing the parameter.
Thus, (\hyperlink{p:ilfp}{InvLFP}) is defined as follows for the optimistic and pessimistic scenario.

{\fontsize{9}{11}
\begin{align}
    & \underbar{\textsc{Inverse Linear Feasibility Problem (optimistic)}} \nonumber \\
    & \left. \begin{aligned}
        &\text{Given:}      && A\in\Q^{m\times n},\ B\in\Q^{m\times k},\ b\in\Q^m,\ C\in\Q^{n\times k}, c\in\Q^n, \\[-0.75ex]
        &                   && \text{a polyhedron } X\subseteq\R^k, \text{ and a set } Y\subseteq\R^n. \\[0.25ex]
        &\text{Question:}   && \text{Does there exist } x\in X \text{ such that} \\[-0.75ex]
        &                   &&  \arg\min\,\{(Cx+c)^\top y \mid Ay \leq Bx+b\}\cap Y\neq \emptyset\;?
    \end{aligned} \right. \tag{InvLFP$_\text{opt}$} \label{p:ilfpo}
\end{align}}
{\fontsize{9}{11}
\begin{align}
    & \underbar{\textsc{Inverse Linear Feasibility Problem (pessimistic)}} \nonumber \\
    & \left. \begin{aligned}
        &\text{Given:}      && A\in\Q^{m\times n},\ B\in\Q^{m\times k},\ b\in\Q^m,\ C\in\Q^{n\times k}, c\in\Q^n, \\[-0.75ex]
        &                   && \text{a polyhedron } X\subseteq\R^k, \text{ and a set } Y\subseteq\R^n. \\[0.25ex]
        &\text{Question:}   && \text{Does there exist } x\in X \text{ such that} \\[-0.75ex]
        &                   &&  \,\emptyset\neq\arg\min\,\{(Cx+c)^\top y \mid Ay \leq Bx+b\}\subseteq Y\;?
    \end{aligned} \right. \tag{InvLFP$_\text{pess}$} \label{p:ilfpp}
\end{align}}
\smallskip

\noindent In the optimistic scenario, we are thus interested in a parameter $x$ such that at least one optimal solution to (\hyperref[p:lpx]{LP$_x$}) lies in $Y$.
In the pessimistic scenario, we are interested in a parameter $x$ such that all optimal solutions to (\hyperref[p:lpx]{LP$_x$}) lie in~$Y$. In both cases, at least one optimal solution must exist.
We will always assume the polyhedron~$X$ to be given as 
$X=\{x\in\R^k \mid Dx\leq d\}$, i.e., that~$D\in\Q^{g\times k}$ and~$d\in\Q^g$ are part of the input of an instance.
As already mentioned, at some points we will instead consider formulations of~(\hyperlink{p:ilfp}{InvLFP}) with an LP in standard form, as we will see that the specific formulation can have a significant influence on the complexity.
If the form of the LP is not specified explicitly in a result, then it applies to both forms. Otherwise the form will be pointed out.

Note that (\hyperlink{p:ilfp}{InvLFP}) could also be defined equivalently without restricting the parameter set.
In fact, for given~$X=\{x\in\R^k \mid Dx\leq d\}$, we could instead allow any parameter~$x\in\R^k$ and, depending on the context, either add~$0\le d-Dx$ as a constraint or add a variable~$\lambda\geq 0$ with the objective function coefficient~$d-Dx$ to the parametric LP. However, this might affect the structure of the target set or the influenced parameters.

\paragraph{Matrix notation.}
We write $[n]\coloneqq \{1,\dots,n\}$ for~$n\in\N$.
For~$A\in\R^{m\times n}$, we denote the $i$-th row vector by $A_{(i)}$ and the $j$-th column vector by $A_{[j]}$.
For index sets~\mbox{$I\subseteq[n]$} and~$J\subseteq[m]$, we denote by $A_{[I]}$ the matrix~$A$ restricted to the columns~$i\in I$ and by $A_{(J)}$ the matrix~$A$ restricted to the rows $j\in J$.
Given parameters $A$, $b$, and~$c$, we will denote the feasible region of the corresponding LP in natural form by $P(A,b) \coloneqq  \{x\in\R^n \mid Ax \leq b\}$ and the feasible region of the corresponding LP in standard form by $P'(A,b) \coloneqq  \{x\in\R_+^n \mid Ax = b\}$, both of which are polyhedra.

A \emph{basis} of $A$ is an index set $\mathcal B\subseteq[n]$ such that~$A_{[\mathcal B]}$ is regular; the associated \emph{(primal) basic solution} for an LP in standard form is~$x\in\R^n$ with~$x_\mathcal B = A_{[\mathcal B]}^{-1} b$ and~$x_{\mathcal N} = 0$, where $\mathcal N\coloneqq [n]\setminus\mathcal B$.
A basis is \emph{feasible} (\emph{optimal}) if the associated basic solution is feasible (resp.\ optimal) for the given LP.
The \emph{dual basic solution} of $\mathcal B$ is~$\mu\in\R^m$ with~$(A^\top\mu)_\mathcal B = c_\mathcal B$.
A basis is \emph{dual feasible} (\emph{dual optimal}) if the dual basic solution is feasible (resp.\ optimal) for the dual LP, which is in natural form.
As we will mostly consider polyhedra with parametrized right-hand sides, we will use the term basic solution independently of the feasibility of the basis, since the latter is only determined once the parameter is fixed. 

\paragraph{Complementary slackness.}
For the proofs of our results, we use some well-known and some less common (unique) optimality conditions for LPs based on complementary slackness.
Consider an LP in natural form
\begin{equation}
    \left.\begin{aligned}
        \min&\quad  c^\top x \\
        \st&\quad    Ax \leq b
    \end{aligned}\qquad\right\} \tag{LP} \label{p:lp}
\end{equation}
and its dual program in standard form
\begin{equation}
    \left.\begin{aligned}
        \max&\quad  b^\top \mu \\
        \st&\quad   A^\top\mu = c,\ \mu \leq 0
    \end{aligned}\qquad \right\} \tag{DLP} \label{p:dlp}
\end{equation}
with $A\in\R^{m\times n}$, $b\in\R^m$, and $c\in\R^n$.
For a primal feasible solution~$x\in\R^n$, a~constraint or index~$j\in[m]$ is called \emph{active} if it is satisfied with equality by~$x$,~i.e., if~$(Ax)_j = b_j$.
We denote the index set of the active constraints for some solution $x$ of \eqref{p:lp} by 
$$\mathcal A(x) \coloneqq  \{j\in[m] \mid (Ax)_j = b_j\}.$$
Similarly, we will refer to the active constraints or indices of an entire face of a polyhedron with the analogous definition.
Vice versa, we refer to the primal face induced by an index set $\mathcal J\subseteq[m]$ by 
$$\face(\mathcal J) \coloneqq  \{x\in\R^n \mid Ax\leq b,\ A_{(\mathcal J)}x=B_\mathcal J\}\;.$$

As first described by Dantzig in the late 1940s, and a special case of the KKT conditions for linear programming, it is well established that a primal feasible solution~$x$ and a dual feasible solution~$\mu$ are each optimal for their respective program if and only if they satisfy the \emph{complementary slackness} condition
\begin{equation}
    (b-Ax)^\top \mu = 0.  \label{eq:complslack}
\end{equation}
With strong duality, this results in well-known necessary and sufficient optimality conditions for both the primal and dual program, which we briefly recall for the convenience of the reader.
\begin{lemma} \label{lem:complementary-slackness-point}
    Consider a primal-dual pair of LPs of the form \eqref{p:lp} and \eqref{p:dlp}.
    If~$x^*$ is feasible for \eqref{p:lp},
    then the following statements are equivalent:
    \begin{itemize}[align=left,labelwidth={1.25em},labelsep=0.25em]
        \item[(i)] $x^*$ is optimal for \eqref{p:lp}.
        \item[(ii)] There exists $\mu$ with $A^\top\mu= c$ and $\mu\leq 0$ such that $x^*$ and $\mu$ satisfy
        \eqref{eq:complslack}.
        \item[(iii)] For every dual optimal solution $\mu^*$, the pair~$(x^*,\mu^*)$ satisfies~\eqref{eq:complslack}.
    \end{itemize}
    Moreover, if $\mu^*$ is feasible for \eqref{p:dlp}, 
    the following statements are equivalent:
    \begin{itemize}[align=left,labelwidth={1.25em},labelsep=0.25em]
        \item[(iv)] $\mu^*$ is optimal for \eqref{p:dlp}.
        \item[(v)] There exists $x$ with $Ax\leq b$ such that $\mu^*$ and $x$ satisfy
        \eqref{eq:complslack}.
        \item[(vi)] For every primal optimal solution $x^*$, the pair~$(x^*,\mu^*)$ satisfies~\eqref{eq:complslack}.
    \end{itemize}
\end{lemma}

Lemma~\ref{lem:complementary-slackness-point} does not provide any sufficient condition for \emph{unique} optimality, which will be important for the analysis of the pessimistic scenario of~\mbox{(\hyperlink{p:ilfp}{InvLFP})}.
To fill this gap, we will establish some equivalent conditions for unique optimality, based on the stronger concept of \emph{strict complementarity}.
\begin{definition}
    A pair $(x,\mu)$ of feasible solutions to \eqref{p:lp} and \eqref{p:dlp} are said to satisfy \emph{strict complementary slackness} if
    \smallskip
    \begin{equation}
        (Ax)_j = b_j \quad\gdw\quad \mu_j\neq 0 \qquad \forall j=1,\dots, m. \label{eq:strictcs}
    \end{equation}
\end{definition}
It is well known that if a primal and dual LP are both feasible, there exists a pair of solutions $(x,\mu)$ that satisfies~\eqref{eq:strictcs}.
Goldman and Tucker~\cite{goldman_tucker1957theory_of_lp} proved this for a different form of LP, but it follows immediately for primal and dual LPs in natural and standard form as well; see \cite{schrijver1998theory}.
From this fact, it is easy to infer a sufficient and necessary condition for the unique optimality of an entire face.

\begin{lemma} \label{lem:uniqueoptface}
    Let $F\neq\emptyset$ be a face of $P(A,b)$ and $\mathcal A(F)$ the set of indices of  active constraints for~$F$.
    Then the following statements are equivalent:
    \begin{itemize}[align=left,labelwidth={1.25em},labelsep=0.25em]
        \item[(i)] $F$ is the set of optimal solutions to \eqref{p:lp}, i.e., $F = \arg\min\{c^\top x \mid Ax\leq b\}$.
        \item[(ii)] There exists a dual feasible solution $\mu$ with
        $\mu_j = 0$ if and only if~$j\not\in \mathcal A(F)$.
    \end{itemize}
\end{lemma}

\begin{proof}
    For the first implication, let $F$ be the set of optimal solutions to \eqref{p:lp}.
    As discussed above, there exists a pair $(x,\mu)$ of primal and dual feasible solutions that satisfies strict complementary slackness.
    Since $x\in F$, we thus have~$\mu_j<0$ for $j\in \mathcal A(F)$.
    On the other hand, for every $j\not\in \mathcal A(F)$, there must exist~$x^{(j)}\in F$ with $(Ax^{(j)})_j < b_j$ by definition of~$\mathcal A(F)$.
    By assumption, this~$x^{(j)}$ is also primal optimal.
    Thus, $(x^{(j)},\mu)$ satisfy complementary slackness as well, and~$\mu_j = 0$ must hold.
    Hence $\mu$ satisfies the condition in~(b).
    
    For the other implication, choose~$\mu$ as in~(b).
    Since~$\mu_j=0$ for~$j\not\in \mathcal A(F)$, the pair~$(x,\mu)$ satisfies complementary slackness for every $x\in F$,
    which implies that~$F\subseteq \arg\min\{c^\top x \mid Ax\leq b\}$.
    Since~$F\neq\emptyset$, this also shows dual optimality of~$\mu$.
    Now let~$x$ be any primal optimal solution.
    This is equivalent to~$x$ being feasible and~$(x,\mu)$ satisfying complementary slackness.
    By assumption on~$\mu$, complementary slackness is only satisfied if~$(Ax)_j=b_j$ for all~$j\in \mathcal A(F)$.
    However, since $F$ is a face of $P(A,b)$, it is uniquely defined by its active constraints. Therefore, $x\in F$ and thus $\arg\min\{c^\top x\mid Ax\leq b\}\subseteq F$.
    
\end{proof}

As a direct consequence of Lemma \ref{lem:uniqueoptface}, we also obtain an equivalent condition for the unique optimality of a single solution.

\begin{corollary}\label{cor:uniqueoptsol}
    Let $x^*$ be a feasible solution to~\eqref{p:lp}.
    Then the following statements are equivalent:
    \begin{itemize}[align=left,labelwidth={1.25em},labelsep=0.25em]
        \item[(i)] $x^*$ is the unique optimal solution to \eqref{p:lp}.
        \item[(ii)] $x^*$ is a vertex of $P(A,b)$ and there exists a dual feasible solution $\mu$ such that~$x^*$ and $\mu$ satisfy strict complementary slackness \eqref{eq:strictcs}.
    \end{itemize}
\end{corollary}

\section{Target solution}\label{sec:ciflp}

We first consider the case where the target set is a singleton, i.e., we ask if there exists a parameter $x\in X$ such that a given \emph{target solution}~$\bar y\in\Q^n$ becomes optimal in the optimistic scenario and uniquely optimal in the pessimistic scenario. 
We will call the resulting special case of~(\hyperlink{p:ilfp}{InvLFP}) the \textsc{Singleton Inverse Linear Feasibility Problem} \hypertarget{p:ilfp-S}{(InvLFP-S)}; it essentially corresponds to the standard case of inverse linear optimization where a single target solution is given.
We begin with a negative result. 

\begin{theorem}[Natural form, RHS] \label{thm:compl-nat}
    Problem \textnormal{(\hyperlink{p:ilfp-S}{InvLFP-S})} is $\npoly$-complete in the optimistic and pessimistic scenario if the LP is given in natural form~\eqref{p:lpx}, even if the objective function is fixed, i.e., if $C=0$.
\end{theorem}
\begin{proof}
    To show $\npoly$-hardness, we first reduce 3-SAT to (\hyperlink{p:ilfp-S}{InvLFP-S$_\text{opt}$}), similarly to the well-known $\npoly$-hardness proof for bilevel optimization \mbox{by Jeroslow~\cite{jeroslow1985bilevel}}.
    So let there be given an instance of 3-SAT with Boolean variables $\mathcal X = \{\xi_1,\dots, \xi_n\}$ and clauses~{$\mathcal C = \{\gamma_1,\dots, \gamma_m\}$}.
    We set $X\coloneqq  [0,1]^n$, where~$x\in X$ is meant to describe a truth assignment for the variables in~$\mathcal X$.
    Further, we define the following LP parametrized by $x$:
    \begin{equation}
        \left. \begin{aligned}
            \min_{y,z} \quad&    -\sum_{i=1}^n y_i - z \\
            \st \quad           &   y_i \leq  x_i,\ y_i \leq 1-x_i                 &&     \forall i=1,\dots, n\\
            &                 z \leq  c_{j1}+c_{j2}+c_{j3}  &&     \forall j=1,\dots, m\\
            &                 z \leq  1 \\
            &                 y \geq 0,\ z\geq 0.
        \end{aligned}\qquad \right\} \label{p:compl-nat-reduction}
    \end{equation}
    In the second set of constraints, the entry $c_{jk}$ represents variable~$x_i$ if the~$k$-th literal in~$\gamma_j$ is~$\xi_i$, while it represents the linear expression~$1-x_i$ if the~$k$-th literal in~$\gamma_j$ is the negation of~$\xi_i$.
    Finally, we define the target solution as~$(\bar y,\bar z)=(0,1)$.

    Now, if there exists a satisfying truth assignment for the instance of 3-SAT, it is easily verified that the target solution is obtained by setting~$x_i=1$ if~$\xi_i$ is true and~$x_i=0$ otherwise. So, conversely, assume that there exists~$x\in X$ such that~$(0,1)$ is optimal for~\eqref{p:compl-nat-reduction}.
    Then $x$ must be binary, as otherwise $y=0$ would not be optimal.
    Thus, $x$ corresponds to a truth assignment for $\mathcal X$, and the latter satisfies all clauses of $\mathcal C$ due to $z=1$.
    This concludes the proof of $\npoly$-hardness for (\hyperlink{p:ilfp-S}{InvLFP-S$_\text{opt}$}).
    The same reduction works for~(\hyperlink{p:ilfp-S}{InvLFP-S$_\text{pess}$}), as it is easy to see that 
    \eqref{p:compl-nat-reduction} has a unique optimal solution for every choice of~$x$.
    
    To prove membership in $\npoly$, we first consider the pessimistic scenario.
    So assume that~$(A,B,b,C,c,X, \bar y)$ is a yes-instance, i.e., there exists $\bar x\in X$ such that~$\bar y$ is the unique optimal solution to $\min\{(Cx+c)^\top y\mid Ay \leq B\bar x + b\}$.
    We choose as certificate the set 
    $$\mathcal A \coloneqq  \{j\in[m] \mid (A\bar y)_j = (B\bar x+b)_j\}$$
    of active indices for~$\bar y$ and parameter~$\bar x$.
    Since $\bar y$ is uniquely optimal for $\bar x$, Corollary~\ref{cor:uniqueoptsol} implies that $\bar y$ is a vertex, i.e., $\rank A_{(\mathcal A)}=n$, and that there exists a dual optimal solution~$\bar\mu$ for parameter $\bar x$ such that~$\bar y$ and~$\bar\mu$ satisfy strict complementary slackness~$(B\bar x+b-A\bar y)^\top\bar\mu = 0$ and~$\bar\mu_{\mathcal A}<0$.
    Now for the verification of the certificate, consider the following LP:
    \refstepcounter{equation} \label{eq:comp-nat-scs}
    \begin{alignat*}{3}
            \max_{x,\mu,\varepsilon} \quad & \varepsilon\qquad \nonumber \\
            \st \quad & 
            (A\bar y)_j = (Bx+b)_j        && \forall j\in\mathcal A \label{eq:comp-nat-scs-a}\tag{\theequation a} \\
            & (A\bar y)_j+\varepsilon \leq (Bx+b)_j  
            && \forall j\not\in\mathcal A \label{eq:comp-nat-scs-b}\tag{\theequation b} \\
            & A^\top\mu = Cx+c,\ \mu\leq 0 \qquad\label{eq:comp-nat-scs-c}\tag{\theequation c} \\
            & \mu_j = 0                     && \forall j\not\in\mathcal A \label{eq:comp-nat-scs-d}\tag{\theequation d} \\
            & \mu_j + \varepsilon \leq 0       && \forall j\in\mathcal A \label{eq:comp-nat-scs-e}\tag{\theequation e} \\
            & x \in X. \label{eq:comp-nat-scs-f}\tag{\theequation f}
    \end{alignat*}
    From our prior observations, it follows that the optimal value of \eqref{eq:comp-nat-scs} is positive, 
    since there exists an $\bar\varepsilon>0$ such that $(\bar x,\bar\mu,\bar\varepsilon)$ satisfies all constraints.
    On the other hand, let~$(x,\mu,\varepsilon)$ be any feasible solution to \eqref{eq:comp-nat-scs} with~$\varepsilon>0$.
    Then~\mbox{\eqref{eq:comp-nat-scs-a}--\eqref{eq:comp-nat-scs-c}} imply that~$\bar y$ is primal feasible and $\mu$ is dual feasible for parameter $x$, while
    \eqref{eq:comp-nat-scs-a} and \eqref{eq:comp-nat-scs-d} imply that $\bar y$ and $\mu$ satisfy complementary slackness and, due to~\eqref{eq:comp-nat-scs-b} and~\eqref{eq:comp-nat-scs-e}, even strict complementary slackness.
    Using Corollary~\ref{cor:uniqueoptsol}, this yields unique optimality of~$\bar y$ for the parameter $x$.
    In summary, the pessimistic instance can be verified with the certificate $\mathcal A$ by showing that $\rank A_{(\mathcal A)}=n$ and that the optimal value of \eqref{eq:comp-nat-scs} is greater than zero.
    
    For the optimistic scenario, let~$\bar y\in\arg\min\{(Cx+c)^\top y\mid Ay \leq B\bar x + b\}$ for some~$\bar x\in X$.
    As before, we choose~$\mathcal A \coloneqq  \{j\in[m] \mid (A\bar y)_j = (Bx+b)_j\}$ as our certificate.
    As~$\bar y$ is optimal for the LP with parameter~$\bar x$, there exists a dual optimal solution $\bar\mu$ such that $\bar y$ and $\bar\mu$ satisfy complementary slackness, i.e.,~$\bar\mu_j=0$ for all~$j\not\in\mathcal A$.
    For the verification, consider the following linear system in~$x$ and~$\mu$:
    \refstepcounter{equation} \label{eq:comp-nat-cs}
    \begin{alignat}{5}
        & A\bar y \leq Bx+b \label{eq:comp-nat-cs-a}\tag{\theequation a} \\
        & (A\bar y)_j = (Bx+b)_j && \qquad \forall j\in\mathcal A \label{eq:comp-nat-cs-b}\tag{\theequation b} \\
        & A^\top\mu = Cx+c,\ \mu\leq 0 \label{eq:comp-nat-cs-c}\tag{\theequation c} \\
        & \mu_j = 0 && \qquad \forall j\not\in\mathcal A \label{eq:comp-nat-cs-d}\tag{\theequation d} \\
        & x\in X\;. \tag{\theequation e} \label{eq:comp-nat-cs-e}
    \end{alignat}
    We already observed that $(\bar x,\bar \mu)$ satisfies \eqref{eq:comp-nat-cs}.
    On the other hand, for any~$(x,\mu)$ that satisfies~\eqref{eq:comp-nat-cs}, it follows that $\bar y$ is optimal for parameter~$x$:
    Inequality~\eqref{eq:comp-nat-cs-a} implies primal feasibility of $\bar y$ and~\eqref{eq:comp-nat-cs-c} implies dual feasibility of $\mu$.
    Together,~\eqref{eq:comp-nat-cs-b} and~\eqref{eq:comp-nat-cs-d} yield complementary slackness.
    In summary, the optimistic instance can be verified with the certificate $\mathcal A$ by showing that \eqref{eq:comp-nat-cs} is feasible.
    
\end{proof}
Theorem \ref{thm:compl-nat} states that the complete case can be hard to decide even if the objective function is certain and only the right-hand side of the LP depends on~$x$.
In fact, the second part of the proof showed that, in this case, the complexity arises from the difficulty of finding the right basis for which $\bar y$ is optimal, out of exponentially many candidates.
In contrast to this, we will show now that the problem becomes easier if only the objective function depends on~$x$.

\begin{theorem}[Natural form, OF]\label{thm:compl-nat-of}
    Problem \textnormal{(\hyperlink{p:ilfp-S}{InvLFP-S})} is tractable in the optimistic and pessimistic scenario if the LP is given in natural form~\eqref{p:lpx} and if the right-hand side is fixed, i.e., if $B=0$.
\end{theorem}
\begin{proof}
    In the optimistic case, $\bar y$ is optimal for $x$ if and only if it is feasible and complementary slackness holds for some dual solution, that is, if and only if the linear system
    $$A\bar y \leq b,
    \quad A^\top\mu = Cx+c, 
    \quad \mu \leq 0,
    \quad (b-A\bar y)^\top \mu = 0,
    \quad x \in X$$
    in $x$ and $\mu$ is feasible. This can be decided in polynomial time.
    
    In the pessimistic case, we know from Corollary~\ref{cor:uniqueoptsol} that~$\bar y$ is the unique optimal solution if and only if it is a feasible vertex and satisfies strict complementary slackness with some dual feasible solution.
    The former property is independent of~$x$ and equivalent to~$A\bar y\leq b$ and $\rank A_{(\mathcal A)} = n$ for
    \mbox{$\mathcal A\coloneqq  \{j\in[m] \mid (A\bar y)_j = b_j \}$}.
    Both can be verified in polynomial time.
    Strict complementary slackness is satisfied if and only if there exists some dual feasible solution~$\mu$ such that~$\bar y$ and~$\mu$ satisfy complementary slackness and $\mu_{\mathcal A}<0$.
    This is equivalent to~the~LP
    \begin{equation*}
        \begin{aligned}
            \max_{x,\mu,\varepsilon} \quad  & \varepsilon \\
            \st \quad & A^\top\mu = Cx+c,\ \mu \leq 0 \\
            & (b-A\bar y)^\top \mu = 0 \\
            & \mu_j+\varepsilon \leq 0 & \forall j\in\mathcal{A} \\
            & x \in X
        \end{aligned}
    \end{equation*}
    having an optimal value greater than zero.
    This LP can be solved efficiently.
\end{proof}
Theorem~\ref{thm:compl-nat} and Theorem~\ref{thm:compl-nat-of} together show that the singleton inverse linear feasibility problem has different complexity status depending on whether the right-hand side or the objective function of the LP is the adjustable parameter.

So far, we have considered~(\hyperlink{p:ilfp-S}{InvLFP-S}) with a parametric LP in natural form. 
Intuitively, one might assume that the complexity does not depend on whether the parametric LP is given in natural or standard form.
However, as already noted by Heuberger~\cite{heuberger2004survey}, the inverse variants of the two forms of LPs cannot be transformed equivalently into each other, since additional slackness variables are added in the transformation of an LP from natural to standard form whose values are not determined by the given natural form solution~$\bar y$.
Instead, we generally end up with a partial inverse linear problem in standard form with a polyhedral target set~$Y = \{\bar y\}\times\R^m$ as discussed in Section~\ref{sec:piflp}.
In fact, the following result proves that (\hyperlink{p:ilfp-S}{InvLFP-S}) can be easier when considering a parametric LP in standard form.

\begin{theorem}[Standard form] \label{thm:comp-stand}
    Problem \textnormal{(\hyperlink{p:ilfp-S}{InvLFP-S})} is tractable in the optimistic and pessimistic scenario if the LP is given in standard form \eqref{p:lpx-stand}.
\end{theorem}
\begin{proof}
    In the optimistic scenario, $\bar y$ is optimal if and only if it is feasible and complementary slackness holds, that is, if and only if the linear system 
    $$ A\bar y = Bx+b, 
    \quad \bar y \geq 0,
    \quad A^\top \mu \leq Cx+c,
    \quad \bar y^\top (Cx+c-A^\top \mu)=0,
    \quad x \in X$$
    in $x$ and~$\mu$ is feasible.
    This can be decided in polynomial time. 
    
    For the pessimistic scenario, according to Corollary \ref{cor:uniqueoptsol}, $\bar y$ is uniquely optimal if and only if it is a feasible vertex and satisfies strict complementary slackness with some dual solution.
    The feasibility of $\bar y$ is equal to $\bar y\geq 0$ and a linear condition in $x$.
    Moreover, to verify that~$\bar y$ is a vertex of the feasible region, we need to check~$\rank A_{[\mathcal I]} = |\mathcal I|$ for the inactive indices~\mbox{$\mathcal I \coloneqq  \{i\in[n] \mid \bar y_i > 0\}$}, which is possible in polynomial time.
    Strict complementary slackness is satisfied if and only if there exists a dual feasible solution~$\mu$ such that~$\bar y$ and~$\mu$ satisfy complementary slackness and~$(Cx+c-A^\top\mu)_{i} >0$ for all~$i\not\in\mathcal I$.
    This is equivalent to the LP
    \begin{equation*}
        \left. \begin{aligned}
            \max_{x,\mu,\varepsilon} \quad  & \varepsilon & \\
            \st \quad & A\bar y = Bx+b \\
            & A^\top\mu \leq Cx+c \\
            & 0=\bar y^\top (Cx+c-A^\top\mu) \\
            & \varepsilon \leq (Cx+c-A^\top\mu)_i &\forall i\not\in{\mathcal I} \\
            & x\in X
        \end{aligned}\qquad \right.
    \end{equation*}
    having an optimal value greater than zero, which can be checked efficiently.
    
\end{proof}
In summary, (\hyperlink{p:ilfp-S}{InvLFP-S}) in standard form is tractable because a complete solution~$\bar y$ already includes all necessary information on how (strict) complementary slackness can be satisfied in order for $\bar y$ to be (uniquely) optimal.
This leads to a possible linearization of the (unique) optimality conditions.
The same approach cannot be used when the LP is in natural form, since the complementary slackness condition would yield a bilinear instead of a linear constraint.

\section{Target basis} \label{sec:targetbasis} \label{subsec:fixed-basis}

Motivated by the positive result of Theorem \ref{thm:comp-stand}, we now consider a modified version of (\hyperlink{p:ilfp-S}{InvLFP-S}) where instead of a specific solution, we fix a basis $\bar{\mathcal B}$ that shall become optimal for some choice of the parameter~$x$: 

\hypertarget{p:ilfp-b}{
{\fontsize{9}{11}\selectfont
\begin{align}
    & \underbar{\textsc{Basis Inverse Linear Feasibility Problem (optimistic)}} \nonumber \\
    & \left. \begin{aligned}
        &\text{Given:}      && A\in\Q^{m\times n},\ B\in\Q^{m\times k},\ b\in\Q^m,\ C\in\Q^{n\times k},\  c\in\Q^n, \\[-0.75ex]
        &                   && \text{a polyhedron } X\subseteq\R^k, \text{ and a basis } \bar{\mathcal B}\text{ of }A. \\[0.25ex]
        &\text{Question:}   && \text{Does there exist } x\in X \text{ such that } \bar{\mathcal B} \text{ is optimal for} \\[-0.75ex]
        &                   && \text{the LP } \min\,\{(Cx+c)^\top y \mid Ay = Bx+b, y\geq 0\}\;?
    \end{aligned} \right. \tag{InvLFP-B$_\text{opt}$} \label{p:ilfp-bo}
\end{align}}
}
\medskip

\noindent The pessimistic version (\hypertarget{p:ilfp-bp}{InvLFP-B$_\text{pess}$}) is defined accordingly by requiring~$\bar{\mathcal B}$ to be the unique optimal basis.
For~(\hyperlink{p:ilfp-b}{InvLFP-B}), we always assume that the parametric LP is given in standard form.

Since a basis always induces a unique basic solution, it might seem plausible that the complexity of~(\hyperlink{p:ilfp-b}{InvLFP-B}) behaves similarly to the complexity of~\mbox{(\hyperlink{p:ilfp-S}{InvLFP-S})}.
Perhaps surprisingly, the following two results prove that there is no direct link between the complexity of the two problems.
Moreover, the complexity even differs between the optimistic and pessimistic scenarios.
\begin{theorem}[Optimistic scenario, RHS]\label{thm:basis-stand-opt}
    Problem \eqref{p:ilfp-bo} is $\npoly$-com\-plete, even if the objective function is fixed.
\end{theorem}
\begin{proof}
    To prove $\npoly$-hardness, we modify the reduction of 3-SAT from the proof of Theorem~\ref{thm:compl-nat}.
    Given an instance $\mathcal X = \{\xi_1,\dots,\xi_n\}$ and $\mathcal C = \{\gamma_1,\dots, \gamma_m\}$ of 3-SAT, we again set~$X\coloneqq [0,1]^n$ and consider a standardized version of \eqref{p:compl-nat-reduction}:
    \begin{equation}
        \left. \begin{aligned}
            \min_{y,z,s} \quad &   -\sum_{i=1}^n y_i - z \\
            \st \quad &  y_i + s_i = x_i    && \forall i=1,\dots, n\\
            & y_i + s_{n+i} = 1- x_i      && \forall i=1,\dots, n \\
            & z + s_{2n+j} = c_{j1}+c_{j2}+c_{j3} &&     \forall j=1,\dots, m\\
            & z + s_{2n+m+1} = 1 \\
            & y \geq 0,\ z\geq 0,\ s\geq 0.
        \end{aligned}\qquad \right\} \label{p:basis-stand-opt-reduction}
    \end{equation}
    As in the proof of Theorem \ref{thm:compl-nat}, it follows that the instance of 3-SAT is satisfiable if and only if there exists a parameter $x\in X$
    such that $(\bar y,\bar z,\bar s)$ with~$\bar y=0$ and~$\bar z=1$ is an optimal solution to \eqref{p:basis-stand-opt-reduction}.
    One easily verifies that the constraint matrix in~\eqref{p:basis-stand-opt-reduction} has full row rank
    and that the variables $z$ and $s_1,\dots,,s_{2n+m}$ induce a basis~$\bar{\mathcal B}$ of~\eqref{p:basis-stand-opt-reduction}. 
    For any parameter $x$, the associated basic solution satisfies~$y={0}$ and~$z=1$.
    Thus, if $\bar{\mathcal B}$ is an optimal and feasible basis, then there exists an optimal solution to~\eqref{thm:basis-stand-opt} with~$y={0}$ and~$z=1$.
    If, on the other hand, there exists such an optimal solution with $\bar y={0}$ and $\bar z=1$, then $\bar{\mathcal B}$ is always a feasible basis that induces this solution.
    In summary, the existence of an optimal solution~$(\bar y,\bar z,\bar s)$ with $\bar y=0$ and~$\bar z=1$ is equivalent to $\bar{\mathcal B}$ being 
    optimal
    for \eqref{p:basis-stand-opt-reduction}.
    This concludes the reduction.
    
    For proving the membership in $\npoly$, let~$(A,B,b,C,c,X,\bar{\mathcal B})$ be a yes-instance of~\mbox{(\hyperlink{p:ilfp-b}{InvFLP-B$_\text{opt}$})}, so that there exists~$\bar x\in X$ such that~$\bar{\mathcal B}$ is an optimal basis for~\eqref{p:lpx-stand}.
    Then the basic solution~$\bar y$ associated with~$\bar{\mathcal B}$ is optimal for parameter~$\bar x$ and there exists a dual optimal solution~$\bar\mu$ such that~$\bar y$ and~$\bar\mu$ satisfy complementary slackness.
    We choose~$\mathcal A \coloneqq  \{i\in[n] \mid \bar y_i=0\}$ as our certificate.
    Then~$\bar x$, $\bar y$, and $\bar \mu$ satisfy the linear system
    \refstepcounter{equation} \label{eq:basis-stand-opt-memb}
    \begin{alignat}{5}
        & Ay = Bx+b \tag{\theequation a},\ y\geq 0 \label{eq:basis-stand-opt-memb-a} \\
        & y_{i} = 0                         && \qquad \forall i\in \mathcal A \tag{\theequation b} \label{eq:basis-stand-opt-memb-b} \\
        & A^\top\mu \leq Cx+c \tag{\theequation c} \label{eq:basis-stand-opt-memb-c} \\
        & (A^\top\mu)_{i} = (Cx+c)_i  && \qquad\forall i\not\in\mathcal A\tag{\theequation d} \label{eq:basis-stand-opt-memb-d} \\
        & x \in X\;, \label{eq:basis-stand-opt-memb-e}\tag{\theequation e}
    \end{alignat}
    since~$i\not\in\mathcal A$ implies $i\in\bar{\mathcal B}$.
    On the other hand, for any subset~$\mathcal A\supseteq[n]\setminus \bar{\mathcal B}$ and for any~$(x,y,\mu)$ satisfying~\eqref{eq:basis-stand-opt-memb}, it follows that $y$ and $\mu$ are primal and dual feasible due to \eqref{eq:basis-stand-opt-memb-a} and \eqref{eq:basis-stand-opt-memb-c}.
    From~\eqref{eq:basis-stand-opt-memb-b} and~\eqref{eq:basis-stand-opt-memb-d} it follows that~$y$ and~$\mu$ satisfy complementary slackness and that~$y$ is the basic solution corresponding to~$\bar{\mathcal B}$.
    Therefore,~$y$ is optimal and, in particular,~$\bar{\mathcal B}$ is an optimal basis for parameter~$x$.
    Thus, the instance can be efficiently verified with the certificate~$\mathcal A$ by showing that~\eqref{eq:basis-stand-opt-memb} is feasible.
    
\end{proof}

Any apparent contradiction 
between Theorem \ref{thm:comp-stand} and Theorem \ref{thm:basis-stand-opt} vanishes when taking a closer look at the proof of Theorem \ref{thm:basis-stand-opt}: Fixing the optimal basic variables for $z$ and $s_1,\dots,s_{2n+m}$ in \eqref{p:basis-stand-opt-reduction}, as done in the reduction, does not yet determine the precise values of the slack variables~$s_1,\dots,s_{2n}$ in the basic feasible solution.
In particular, it is not clear which of the variables~$s_{1},\dots,s_{2n}$ are positive and which are zero.
In fact, if~$x$ is binary, as necessary for~$\bar{\mathcal B}$ to be optimal, then exactly~$n$ of the variables~$s_1,\dots,s_n$ and~$s_{n+1}\dots,s_{2n}$ are zero and could be arbitrarily replaced by the respective variables~$y_i$ to form other feasible bases inducing the same basic solution.
So if~$\bar{\mathcal B}$ is optimal, there exist at least~$2^n$ optimal bases.
Moreover, not all those bases need to be dual feasible, and guessing the dual optimal basis is equivalent to guessing the truth assignment to the Boolean variables.
Thus, in order to check if $\mathcal B$ can be optimal for some parameter~$x\in X$, one needs to check for $2^n$ potential dual optimal bases if they satisfy complementary slackness together with 
$\bar{\mathcal B}$.
On the other hand, if a specific target solution is given, the zero and non-zero values of the complete solution already uniquely determine which variables need to be part of a dual optimal basis, making the problem tractable.

\begin{remark}[Optimistic scenario, OF]\label{rem:basis-stand-opt-of}
    It is easy to verify that~\eqref{p:ilfp-bo} becomes tractable if the right-hand side is fixed and only the objective function is parametrized by~$x$: For a given basis, its basic solution is unique and independent of~$x$. Thus, similar to Theorem~\ref{thm:compl-nat-of}, feasibility of the basis can be checked independently of~$x$ and the complementary slackness condition for the basic solution can be linearized.
    We omit a detailed proof.
\end{remark}

We next prove that~(\hyperlink{p:ilfp-b}{InvLFP-B}) is tractable in the pessimistic scenario.
This result is interesting, since in bilevel optimization it is more common for the pessimistic scenario to be harder than the optimistic scenario.
Moreover, this is one of the few cases of (\hyperlink{p:ilfp}{InvLFP}) we considered for which the complexities of the optimistic and pessimistic scenario differ.

\begin{theorem}[Pessimistic scenario] \label{thm:basis-stand-pess}
    Problem \textnormal{(\hyperlink{p:ilfp-b}{InvLFP-B$_\text{pess}$})} is tractable.
\end{theorem}
\begin{proof}
    We explain how a given instance $(A,B,b,C,c,X,\bar{\mathcal B})$ can be decided in polynomial time.
    If $\bar{\mathcal B}$ is a feasible basis for the parameter $x\in X$, the system
    \begin{equation}
        A y = Bx+b,
        \quad y \geq 0,
        \quad y_{\bar{\mathcal N}} = 0
        \label{eq:basis-stand-pess-y}
    \end{equation}
    is feasible, where~$\bar{\mathcal N}=[n]\setminus\bar{\mathcal B}$ is the non-basis corresponding to~$\bar{\mathcal B}$.
    Let~$\bar y$ be the unique solution to \eqref{eq:basis-stand-pess-y} for the given~$x\in X$. We first need to address the difficulty that different bases might induce the same basic solution.
    For~$\bar{\mathcal B}$ to be uniquely optimal, there cannot exist any index~$i\in\bar{\mathcal B}$ with $\bar y_i = 0$ and some non-basic index~$j\in\bar{\mathcal N}$ such that $(\bar{\mathcal B}\setminus\{i\})\cup\{j\}$ is again a basis of~$A$.
    We define
    \begin{equation*}
        {\bar{\mathcal B}}^{0} \coloneqq  \left\{ i\in\bar{\mathcal B} \mid \exists j\in\bar{\mathcal N}\colon \left(\bar{\mathcal B}\setminus\{i\}\right)\cup\{j\} \text{ is a basis of } A \right\}
    \end{equation*}
    and refer to $\bar{\mathcal B}^{0}$ as the \emph{non-essential part} of $\bar{\mathcal B}$ and to $\bar{\mathcal B}\setminus\bar{\mathcal B}^{0}$ as the \emph{essential part} of $\bar{\mathcal B}$.
    By considering the determinant of the basis matrix, it holds
    \begin{equation*}
        \bar{\mathcal B}^{0} = \{i\in\bar{\mathcal B} \mid \exists j\in\bar{\mathcal N} : \alpha^{(j)}_i \neq 0\}
    \end{equation*}
    where~$\alpha^{(j)}\coloneqq A_{[\bar{\mathcal B}]}^{-1} A_{[j]}$.
    Clearly~$\bar{\mathcal B}^0$ consists precisely of the indices that must be part of every basis of $A$, implying that the essential part is equal for all bases.
    Thus, the basis~$\bar{\mathcal B}$ is the unique basis inducing~$\bar y$ if and only if $\bar y_i > 0$ for all $i\in\bar{\mathcal B}^{0}$.
    Further, we need to prevent the existence of other optimal basic solutions.
    To this end, we set
    \begin{equation*}
        \bar{\mathcal N}^+ \coloneqq \{j\in\bar{\mathcal N} \mid \exists i\in\bar{\mathcal B}^0: \alpha^{(j)}_i>0\}\;.
    \end{equation*}
    We claim that for any~$j'\in\bar{\mathcal N}^+$, there exists a basis~$\mathcal B'$ resulting from pivoting~$j'$ with some index in~$\bar{\mathcal B}^0$ such that~$\mathcal B'$ is feasible if~$\bar{\mathcal B}$ was feasible.
    For this, select
    \begin{equation*}
        i'\in\arg\min\Big\{\tfrac{\bar y_i}{\alpha^{(j')}_i} \mid i\in\bar{\mathcal B}^0, \alpha^{(j')}_i>0\Big\}
    \end{equation*}
    as the leaving index.
    Since~$A_{[j']} = A_{[\bar{\mathcal B}]}\alpha^{(j)}$, the basic solution~$y'$ of~$\mathcal B'$ satisfies~$y'_{j'} = \tfrac{\bar y_{i'}}{\alpha^{(j')}_{i'}}$ and for the remaining basic indices~$i\in\bar{\mathcal B} \setminus \{i'\}$
    \begin{align*}
        y'_i = \bar y_i - y'_{j'}\alpha^{(j')}_i
        &\geq 
        \begin{cases}
            \bar y_i \quad &\text{if } \alpha^{(j')}_i\leq 0\;, \\
            \bar y_i - \tfrac{\bar y_{i'}}{\alpha^{(j')}_{i'}} \alpha^{(j')}_{i} \quad &\text{if } \alpha^{(j')}_i>0\;.
        \end{cases}
    \end{align*}
    Therefore, if~$\bar y\geq 0$, then~$y'\geq 0$, so~$\mathcal B'$ is also feasible.
    Now given $\bar{\mathcal B}$, we can determine $\bar{\mathcal B}^{0}$ and~$\bar{\mathcal N}^+$ in polynomial time.
    Altogether, we claim that $\bar{\mathcal B}$ is uniquely optimal for some feasible parameter $x$ if and only if the following linear system has a positive objective value:
    \refstepcounter{equation}
    \label{eq:basis-stand-pess-opt}
    \begin{alignat*}{3}
        \max_{x,y,\mu,\varepsilon} \quad & \varepsilon\\
        \st \quad & Ay = Bx+b,\ y\geq 0 \tag{\theequation a} \label{eq:basis-stand-pess-opt-a} \\
        & y_i = 0   && \forall i\in\bar{\mathcal N} \tag{\theequation b} \label{eq:basis-stand-pess-opt-b} \\
        & y_i \ge \varepsilon \qquad                      && \forall i\in\bar{\mathcal B}^{0} \tag{\theequation c} \label{eq:basis-stand-pess-opt-c} \\
        & A^\top \mu \leq Cx+c   \tag{\theequation d} \label{eq:basis-stand-pess-opt-d} \\
        & (A^\top\mu)_i = (Cx+c)_i\qquad    && \forall i\in\bar{\mathcal B}^{0} \tag{\theequation e} \label{eq:basis-stand-pess-opt-e} \\
        & (A^\top\mu)_i+\varepsilon \le (Cx+c)_i\quad && \forall i\in\bar{\mathcal N}^+ \tag{\theequation f} \label{eq:basis-stand-pess-opt-f}\\
        & x\in X\tag{\theequation g} \label{eq:basis-stand-pess-opt-g}\;.
    \end{alignat*}
    We need to prove both directions of the equivalence.
    If $\bar{\mathcal B}$ is uniquely optimal for some parameter~$x$, then we already argued that \eqref{eq:basis-stand-pess-opt-a}--\eqref{eq:basis-stand-pess-opt-c} need to be satisfied such that~$\bar{\mathcal B}$ is feasible and the unique basis to produce its basic solution.
    Furthermore, for the optimality, complementary slackness~\eqref{eq:basis-stand-pess-opt-d}--\eqref{eq:basis-stand-pess-opt-e} must be satisfied.
    The dual optimal solution~$\mu$ must also satisfy~\eqref{eq:basis-stand-pess-opt-f}, as otherwise we argued that there exists another basic feasible solution which would also satisfy complementary slackness with~$\mu$ and therefore be optimal.
    The constraint~\eqref{eq:basis-stand-pess-opt-g} is trivial.
    For the other direction, let~$(x,y,\mu,\varepsilon)$ be a feasible solution to~\eqref{eq:basis-stand-pess-opt} with~\mbox{$\varepsilon>0$}.
    Then~\eqref{eq:basis-stand-pess-opt-a}--\eqref{eq:basis-stand-pess-opt-c} imply that $\bar{\mathcal B}$ is feasible and the unique basis to induce~$y$.
    It remains to derive the optimality of~$\bar{\mathcal B}$ and the non-existence of another optimal basic solution.
    The tuple~$(y,\mu)$ satisfies complementary slackness in all but the essential indices~$i\in\bar{\mathcal B}\setminus\bar{\mathcal B}^0$.
    However, since~$A_{[i]}$ is independent from~$A_{[[n]\setminus\{i\}]}$ for every~$i\in\bar{\mathcal B}\setminus\bar{\mathcal B}^0$, the set~$\{z\in\R^m \mid A_{([n]\setminus\{i\})}^\top z=0,\ A_{(i)}^\top z = \delta\}$ is non-empty for every~$i\in\bar{\mathcal B}\setminus\bar{\mathcal B}^0$ and~$\delta\in\R$.
    Thus,~$\mu$ can be extended to a dual feasible solution~$\mu'$ that satisfies, additionally to~\eqref{eq:basis-stand-pess-opt-d}--\eqref{eq:basis-stand-pess-opt-f}, also~$(A^\top\mu)_i=(Cx+c)_i$ for all~$i\in\bar{\mathcal B}\setminus\bar{\mathcal B}^0$.
    Then~$(y,\mu')$ satisfies complementary slackness everywhere and~$y$ is optimal.
    Let~$z$ be another primal optimal solution.
    Then it holds
    \begin{align*}
        z_{\bar{\mathcal B}} - y_{\bar{\mathcal B}}
        &= A_{[\mathcal{\bar B}]}^{-1}(Bx+b-A_{[\mathcal{\bar N}]}z_{\bar{\mathcal N}}) - A_{[\bar{\mathcal B}]}^{-1}(Bx+b) \\
        &= -A_{[\bar{\mathcal B}]}^{-1} A_{[\bar{\mathcal N}]} z_{\bar{\mathcal N}} \\
        &= -\sum_{i\in\bar{\mathcal N}\setminus\bar{\mathcal N}^+} A_{[\bar{\mathcal B}]}^{-1} A_{[i]} z_i \\
        &= -\sum_{i\in\bar{\mathcal N}\setminus\bar{\mathcal N}^+} \alpha^{(i)}z_i \geq 0\;.
    \end{align*}
    The third equation holds since, due to~\eqref{eq:basis-stand-pess-opt-f}, any optimal solution has to satisfy~$z_{\bar{\mathcal N}^+}=0$.
    The last inequality follows because~$\alpha^{(j)}\leq 0$ for all~$j\in\bar{\mathcal N}\setminus\bar{\mathcal N}^+$ by the definition of~$\bar{\mathcal N}^+$.
    The inequality, however, implies that~$z - y\geq 0$ is a feasible direction from~$y$, so the entire ray from~$y$ in the direction~$z-y$ is optimal.
    But then~$z$ cannot be a vertex of the feasible polyhedron, so~$y$ is the only optimal basic solution.
    By~\eqref{eq:basis-stand-pess-opt-g}, the parameter~$x$ is feasible.
    
    
\end{proof}

The proof of Theorem~\ref{thm:basis-stand-pess} demonstrates why, different from the optimistic scenario, the pessimistic scenario is tractable: A yes-instance is now characterized by the fact that~$\bar{\mathcal B}$ is the unique optimal basis, implying that $\bar{\mathcal B}$ is (essentially) non-degenerate, whereas degeneracy was the main issue that led to the hardness of the optimistic scenario.
Moreover, if~$\bar{\mathcal B}$ is the unique primal optimal basis, it must also be dual optimal.

\section{Polyhedral target set} \label{sec:piflp}

We now consider~(\hyperlink{p:ilfp}{InvLFP}) with target sets~$Y$ containing more than one element.
First, note that the case of finite~$Y$ can easily be reduced to the complete case discussed in Section \ref{sec:ciflp} and vice versa, so that the corresponding complexity results carry over to arbitrary target sets given explicitly as a list in the input.

We will thus focus on target \emph{polyhedra} in the remainder of this section. So assume that~$Y = \{y\in\R^n \mid Sy\leq t\}$, where~$S\in\Q^{r\times n}$ and~$t\in\Q^r$ are now part of the input of~(\hyperlink{p:ilfp}{InvLFP}).
Then the $\npoly$-hardness of~(\hyperlink{p:ilfp}{InvLFP}) follows from the results in Section~\ref{sec:ciflp}, but also from the $\npoly$-hardness of bilevel linear programming without coupling constraints. Indeed, the decision variant of the latter can be modeled as (\hyperlink{p:ilfp}{InvLFP})
by defining the target polyhedron as the set of all~$y$ with a leader's objective value below a certain threshold.
In fact, we will show in this section that $\npoly$-hardness already holds for much simpler polyhedra.

However, we begin with a general positive result.
\begin{theorem}[Target polyhedron] \label{thm:polyh-np}
    If~$Y$ is a polyhedron, then \textnormal{(\hyperlink{p:ilfp}{InvLFP})} belongs to $\npoly$ both in the optimistic and in the pessimistic scenario.
\end{theorem}
\begin{proof}
    First, consider a yes-instance $(A,B,b,C,c,X,Y)$ of \eqref{p:ilfpo}, so that there exists $\bar x\in X$ with $Y\cap \bar F\neq \emptyset$, where
    $$\bar F\coloneqq \arg\min\{(C\bar x+c)^\top y \mid Ay\leq B\bar x+b\}\;.$$
    In particular, $\bar F$ is a non-empty face of $P(A,B\bar x+b)$. Our certificate will consist of the set of active indices
    $$\mathcal A \coloneqq  \{j\in[m] \mid (Ay)_j=(B\bar x+b)_j\ \forall y\in \bar F\}$$
    for~$\bar F$ and the parameter~$\bar x$.
    From Lemma~\ref{lem:uniqueoptface}, we know that there exists a dual feasible solution $\bar\mu$ with $\bar\mu_j = 0$ if and only if $j\not\in \mathcal A$.
    As $Y\cap \bar F\neq \emptyset$, there exists some $\bar y\in Y$ with $A\bar y\leq B\bar x+b$ and $(A\bar y)_j = (B\bar x+b)_j$ for all $j\in \mathcal A$.
    Altogether, the following linear system is satisfied by $(\bar x,\bar y,\bar \mu)$:
    \refstepcounter{equation} \label{eq:polyh-np-opt}
    \begin{alignat}{5}
        & Ay \leq Bx+b \tag{\theequation a} \label{eq:polyh-np-opt-a} \\
        & (Ay)_j = (Bx+b)_j \quad && \forall j\in \mathcal A \tag{\theequation b} \label{eq:polyh-np-opt-b} \\
        & A^\top\mu = Cx+c \tag{\theequation c} \label{eq:polyh-np-opt-c} \\
        & \mu \leq 0 \tag{\theequation d} \label{eq:polyh-np-opt-d} \\
        & \mu_j = 0 && \forall j\in[m]\setminus \mathcal A \tag{\theequation e} \label{eq:polyh-np-opt-e} \\
        & Sy \leq t \tag{\theequation f} \label{eq:polyh-np-opt-f} \\
        & x \in X. \tag{\theequation g} \label{eq:polyh-np-opt-g}
    \end{alignat}
    On the other hand, for any triple $(x,y,\mu)$ that solves \eqref{eq:polyh-np-opt}, it follows that $x$ is a solution to the given instance of \eqref{p:ilfpo}. Indeed, by Lemma~\ref{lem:uniqueoptface}, the constraints~\eqref{eq:polyh-np-opt-c}--\eqref{eq:polyh-np-opt-e} imply that the face induced by the constraints $j\in \mathcal A$ is a subset of the set of optimal solutions.
    Due to~\eqref{eq:polyh-np-opt-a}--\eqref{eq:polyh-np-opt-b}, the face induced by~$\mathcal A$ contains~$y$, and~\eqref{eq:polyh-np-opt-f} yields $y\in Y$.
    The last constraint \eqref{eq:polyh-np-opt-g} ensures that~$x$ is a feasible parameter for the instance.

    Next, let $(A,B,b,C,c,X,Y)$ be a yes-instance of~\eqref{p:ilfpp}, so that there exists $\bar x\in X$ with $\emptyset \neq \bar F \subseteq Y$, where~$\bar F$ is defined as above.
    The set of active indices~$\mathcal A$ will be the first part of our certificate.
    The verification can be performed similarly to the optimistic scenario, with the exception that we now must verify that the face of $P(A,Bx+b)$ induced by $\mathcal A$ is a subset of~$Y$.
    We denote this face, which depends on~$x$, by
    $$\face (\mathcal A,x) \coloneqq  \{y\in \R^n \mid Ay\leq Bx+b,\ -(Ay)_{\mathcal A} \le -(Bx+b)_{\mathcal A}\}\;.$$
    Since~$\bar F\subseteq Y$, every constraint defining~$Y$ is valid for~$\bar F = \face (\mathcal A, \bar x)$.
    Thus, for every~$h\in[r]$, there exist slack vectors~$s^{(h,1)}\in\R_+^m$ and $s^{(h,2)}\in\R_+^{\mathcal A}$ such that
    \begin{equation}
        \left.\begin{aligned}
        & S_{(h)} = A^\top s^{(h,1)} - A_{(\mathcal A)}^\top s^{(h,2)}\\
        & t_h \geq (B\bar x+b)^\top s^{(h,1)} - (B\bar x+b)_{\mathcal A} ^\top s^{(h,2)}\;.
        \end{aligned}\qquad \right\} \label{eq:conic_comb2}
    \end{equation}
    For each $h\in[r]$, we set $${\mathcal S}_+^{(h,1)}\coloneqq  \{j\in[m] \mid s^{(h,1)}_j > 0\},\quad {\mathcal S}_+^{(h,2)} \coloneqq  \{j\in \mathcal A \mid s^{(h,2)}_j > 0\}\;.$$
    We may assume
    ${\mathcal S}_+^{(h,1)}\cap {\mathcal S}_+^{(h,2)}=\emptyset$ and that the vectors $A_{(j)}$, $j\in {\mathcal S}_+^{(h,1)}\cup {\mathcal S}_+^{(h,2)}$, are linearly independent for every $h\in[r]$. These~$2r$ index sets will be the second part of the certificate.
    With this assumption, however, for every $h\in[r]$ the pair~$({s}^{(h,1)},{s}^{(h,2)})$ is the unique solution to the system
    \begin{equation}
        \left.\begin{aligned}
            & S_{(h)}             = A^\top {s}^{(h,1)} - A_{(\mathcal A)}^\top {s}^{(h,2)} \\
            & {s}^{(h,1)}_j     = 0\qquad      &&\forall j\not\in {\mathcal S}_+^{(h,1)} \\
            & {s}^{(h,2)}_j     = 0\qquad      &&\forall j\not\in {\mathcal S}_+^{(h,2)} \\
            & {s}^{(h,1)}       \in \R_+^m,\ {s}^{(h,2)}       \in \R_+^{\mathcal A}
        \end{aligned}\qquad \right\} \label{eq:polyh-np-pess-pos-comb}
    \end{equation}
    and this unique solution also satisfies~\eqref{eq:conic_comb2}.
    For the actual verification, we consider the following LP:
    \pagebreak[3]
    \refstepcounter{equation} \label{eq:polyh-np-pess}
    \begin{alignat*}{5}
        \max_{x,y,\mu,\varepsilon} \quad & \varepsilon \\
        \st \quad & (Ay)_j = (Bx+b)_j                 && \forall j\in \mathcal A \tag{\theequation a} \label{eq:polyh-np-pess-a} \\
        & (Ay)_j+\varepsilon \leq (Bx+b)_j                 && \forall j\not\in \mathcal A \tag{\theequation b} \label{eq:polyh-np-pess-b} \\
        & A^\top\mu = Cx+c \tag{\theequation c} \label{eq:polyh-np-pess-c} \\
        & \mu_j = 0                         && \forall j\not\in \mathcal A \tag{\theequation d} \label{eq:polyh-np-pess-d} \\
        &  \mu_j+\varepsilon\leq 0 && \forall j\in \mathcal A \tag{\theequation e} \label{eq:polyh-np-pess-e} \\
        & t_h \geq (Bx+b)^\top {s}^{(h,1)} - (Bx+b)_{\mathcal A}^\top {s}^{(h,2)}\quad  && \forall h\in[r] \tag{\theequation f} \label{eq:polyh-np-pess-f} \\
        & x \in X. \tag{\theequation g} \label{eq:polyh-np-pess-g}
    \end{alignat*}
    As in the optimistic scenario, for a yes-instance there exists a primal optimal solution $\bar y\in\bar F$ and a dual optimal solution $\bar\mu$ for parameter $\bar x$ that satisfy strict complementary slackness.
    Thus, there exists an $\bar\varepsilon>0$ such that $(\bar x, \bar y, \bar \mu,\bar\varepsilon)$ is a feasible solution to \eqref{eq:polyh-np-pess}.
    On the other hand, we claim that if \eqref{eq:polyh-np-pess} has a feasible solution~$(x,y,\mu,\varepsilon)$ with~$\varepsilon>0$, then~$x$ solves the instance of \eqref{p:ilfpp}.
    Indeed, constraints \eqref{eq:polyh-np-pess-a}--\eqref{eq:polyh-np-pess-b} imply that the face~$\face(\mathcal A,x)$ of $P(A,Bx+b)$ induced by $\mathcal A$ is non-empty, and due to \eqref{eq:polyh-np-pess-c}--\eqref{eq:polyh-np-pess-e} and Lemma~\ref{lem:uniqueoptface}, this face is the unique set of optimal solutions for the LP with parameter $x$.
    Finally,~\eqref{eq:polyh-np-pess-f} yields $\face(\mathcal A,x)\subseteq Y$, showing that $x$ solves the instance of \eqref{p:ilfpp}.

    The verification of the certificate thus consists of 
    computing the vectors~$s^{(h,1)}$ and~$s^{(h,2)}$ by solving the system~\eqref{eq:polyh-np-pess-pos-comb} and then checking whether the LP~\eqref{eq:polyh-np-pess} has a positive optimal value.
    
\end{proof}

Note that the result of Theorem \ref{thm:polyh-np} holds independently of the form of the underlying parametric LP, because the polynomial transformation of (\hyperlink{p:ilfp}{InvLFP}) in standard form into an equivalent instance of (\hyperlink{p:ilfp}{InvLFP}) in natural form preserves the polyhedral property of the target set.

It already follows from Theorem~\ref{thm:compl-nat} that  (\hyperlink{p:ilfp}{InvLFP}) with a polyhedral target set~$Y$ is $\npoly$-hard in general, even in the strictest possible case~$Y=\{\bar y\}$. This raises the question of how hard (\hyperlink{p:ilfp}{InvLFP}) can become in the opposite case, i.e., when $Y=\R^n$.
In other words, we ask whether there exists a parameter~$x\in X$ such that the LP has any optimal solution at all.
However, the latter question is easily answered by checking the feasibility of the following linear system
$$Ay \leq Bx+b,\quad A^\top \mu = Cx+c,\quad \mu \geq 0, \quad x \in X$$
which models primal and dual feasibility.
This in turn raises the question of how far we can restrict~$Y$ such that (\hyperlink{p:ilfp}{InvLFP}) remains tractable.
Unfortunately, the following two theorems show that even restricting only a single variable may already result in an $\npoly$-hard problem.

\begin{theorem}[One target variable, RHS] \label{thm:one-variable-nat-rhs}
    Problem \textnormal{(\hyperlink{p:ilfp}{InvLFP})} is $\npoly$-complete both in the optimistic and in the pessimistic scenario if $Y = \{y\in\R^n \mid y_i= 0\}$ for some $i\in[n]$, even if the objective function is fixed.
\end{theorem}
\begin{proof}
    Membership in $\npoly$ follows from Theorem~\ref{thm:polyh-np}.
    To prove $\npoly$-hardness of the optimistic scenario, we first describe a Karp-reduction of 3-SAT to \eqref{p:ilfpo} which is very similar to the one used in the proof of Theorem \ref{thm:compl-nat}.
    Given an instance of 3-SAT with $n$ Boolean variables, we set~$X \coloneqq  [0,1]^n$ and consider the parametrized~LP
    \begin{equation}
        \left. \begin{aligned}
            \min_{y,z,t} \quad &  -\sum_{i=1}^n y_i - z\qquad\qquad \\
            \st \quad & y_i \leq x_i,\; y_i \leq 1-x_i &&     \forall i=1,\dots, n \\
            & z \leq c_{j1}+c_{j2}+c_{j3}  &&     \forall j=1,\dots, m\\
            & z \leq 1 \\
            & -\sum_{i=1}^n y_i + z + t = 1 \\
            & y \geq  0,\ z\geq  0,\ t \geq 0\;.
        \end{aligned}\qquad \right\} \label{p:one-variable-nat-reduction}
    \end{equation}
    Note that~\eqref{p:one-variable-nat-reduction} agrees with~\eqref{p:compl-nat-reduction}, with the addition of variable~$t\geq 0$ and the second last constraint (which corresponds to two inequalities in natural form). 
    This, however, does not affect the feasibility or optimal value of the program, since the additional constraint can always be satisfied by setting~\mbox{$t = \sum_{i=1}^n y_i -z + 1\geq 0$} and $t$ is not part of the objective function.
    
    Analogously to the proof of Theorem~\ref{thm:compl-nat}, it follows that the instance of 3-SAT is satisfiable if and only if there exists a parameter $x\in X$ such that~$(y^*,z^*,t^*)$ with $y^*=0$ and $z^* =1$ is an optimal solution to \eqref{p:one-variable-nat-reduction}, which implies~$t^* = 0$.
    On the other hand, if $t^*=0$, the second last constraint is only satisfied if~$y^*=0$ and~$z^*=1$.
    In summary, setting~$Y\coloneqq  \{(y,z,t)\in\R^n\times\R\times\R \mid t= 0\}$, the instance of 3-SAT is satisfiable if and only if there exists a parameter $x\in X$ such that there exists an optimal solution $(y^*,z^*,t^*) \in Y$ to \eqref{p:one-variable-nat-reduction}.
    This finishes the reduction and proves the $\npoly$-hardness of \eqref{p:ilfpo}.
    
    The same reduction also proves $\npoly$-hardness of \eqref{p:ilfpp}, as one quickly verifies that~$(y^*,z^*,t^*)=(0,1,0)$ is the unique optimal solution to \eqref{p:one-variable-nat-reduction} if it is an optimal solution at all.
    
\end{proof}

So far, all $\npoly$-hardness results for (\hyperlink{p:ilfp}{InvLFP}) applied to cases where the parameter~$x$ affected the right-hand side of the LP.
Moreover, for (\hyperlink{p:ilfp-S}{InvLFP-S}), we even proved that the problem is tractable if the right-hand side is fixed.
Nevertheless, the next result shows that (\hyperlink{p:ilfp}{InvLFP}) can also be $\npoly$-hard even if only the objective function depends on the parameter~$x$.
\begin{theorem}[One target variable, OF] \label{thm:one-variable-nat-of}
    Problem \textnormal{(\hyperlink{p:ilfp}{InvLFP})} is $\npoly$-complete both in the optimistic and in the pessimistic scenario if $Y = \{y\in\R^n \mid y_i = 0\}$ for some~$i\in[n]$, even if the right-hand side is fixed.
\end{theorem}
\begin{proof}
    Membership in $\npoly$ follows from Theorem~\ref{thm:polyh-np}.
    To prove $\npoly$-hardness, we again reduce 3-SAT to (\hyperlink{p:ilfp}{InvLFP}). 
    However, in order to deal with the fixed right-hand side, we now consider a parametric LP that is related to the dual program of~\eqref{p:compl-nat-reduction}.
    Using the same notation as before, we now set
    \begin{align*}
        X \coloneqq  \{ x\in[0,1]^n \mid\ &c_{j1}+c_{j2}+c_{j3} \geq 1\quad \forall j=1,\dots, m\}
    \end{align*}
    and consider the parametric LP
    \begin{equation}
        \left. \begin{aligned}
            \min_{\mu, \nu, z} \quad & -x^\top\mu-({1}_n-x)^\top\nu - \big( n - \tfrac 14 \big) z \\
            \st \quad & \mu_i + \nu_i + z \leq 1\quad\qquad \forall i=1,\dots, n \\
            & \mu\geq 0,\ \nu \geq 0,\ z \geq 0
        \end{aligned}\qquad \right\} \label{p:one-variable-nat-of-reduction}
    \end{equation} 
    as well as the target set $Y \coloneqq  \{(\mu,\nu,z)\in\R^n\times\R^n\times\R \mid z=0\}$.
    We claim that the given instance of 3-SAT is satisfiable if and only if the constructed instance is a yes-instance of~\eqref{p:ilfpo} if and only if it is a yes-instance of~\eqref{p:ilfpp}.
    Since yes-instances of~\eqref{p:ilfpp} are yes-instances of~\eqref{p:ilfpo}, it suffices to show two implications:
    First, we show that satisfiability of the 3-SAT instance implies that the constructed instance is a yes-instance of~\eqref{p:ilfpp}. Then we show that if the instance is a yes-instance of~\eqref{p:ilfpo}, this implies satisfiability of the 3-SAT instance.

    So assume that the 3-SAT instance is satisfiable and define~$x$ through a satisfying truth assignment; then clearly~$x\in X$.
    It is easy to verify that~\mbox{$(x,{1}_n-x,0)$} is the unique optimal solution to~\eqref{p:one-variable-nat-of-reduction} for this~$x$, showing that we have indeed constructed a yes-instance of~\eqref{p:ilfpp}.
    
    Now assume that we have a yes-instance of~\eqref{p:ilfpo}, i.e., there exists a parameter~$x\in X$ such that some~$(\mu,\nu,0)\in Y$ is optimal for~\eqref{p:one-variable-nat-of-reduction}.
    We claim that the truth assignment $T\colon\mathcal X\to\{0,1\}$ with
    \begin{equation*}
        T(\xi_i) = \begin{cases}
                        ~~ 1\quad  & \text{if } x_i \geq \frac 34 \\[0.75ex]
                        ~~ 0\quad  & \text{if } x_i \leq \frac 14
                    \end{cases}
    \end{equation*}
    is well-defined and satisfies the instance of 3-SAT.
    For the first claim, assume on contrary that there exists an $\hat\imath\in[n]$ with $\max\{x_{\hat\imath}, 1- x_{\hat\imath} \} < \frac 34$.
    Then for any feasible solution $(\mu,\nu,0)\in Y$, the objective function value of~\eqref{p:one-variable-nat-of-reduction} is
    \[
        -x^\top\mu-({1}_n-x)^\top \nu \geq -\sum_{i=1}^n \max\{x_i, 1- x_i\} > \tfrac 14-n\;.
    \]
    However, this is a contradiction to the fact that $(0,0,1)$ is always feasible with objective value $\tfrac 14-n$ and that we consider a yes-instance of~\eqref{p:ilfpo}.
    Now if $\max\{x_i,1-x_i\}\geq \frac 34$, then $\min\{x_i,1-x_i\}\leq \frac 14$, since $x_i\in[0,1]$.
    Moreover,~$x\in X$ implies $c_{j1}+c_{j2}+c_{j3}\geq 1$ for every $j\in[m]$, so that at least one of the $c_{j1},c_{j2},c_{j3}$ must be at least~$\frac 34$.
    In other words, $T$ satisfies the given instance of 3-SAT.
    
\end{proof}
The last two results show that the complexity of (\hyperlink{p:ilfp}{InvLFP}) can change from tractable to $\npoly$-hard by fixing a single dimension of the optimization variable~$y$.
This applies regardless of whether the underlying parametric LP is given in natural or standard form, since both~\eqref{p:one-variable-nat-reduction} and~\eqref{p:one-variable-nat-of-reduction} can be transformed into the other form without compromising the proof.

To conclude this section, we now consider the most common case of partial inverse optimization, where the partial solution consists of a subset of solution variables that are fixed in the target set.
The previous two theorems as well as the $\npoly$-hardness of (\hyperlink{p:ilfp-S}{InvLFP-S}) in natural form (see Theorem \ref{thm:compl-nat}) imply that this subproblem of (\hyperlink{p:ilfp}{InvLFP}) will generally be $\npoly$-hard in natural form, for any number of fixed variables.
For standard form, however, we recall that the problem is tractable if all variables are fixed; see Theorem \ref{thm:comp-stand}.
The following result shows that the complexity does in fact not increase drastically if the number of free variables is increased gradually.
Instead, we are able to prove fixed parameter tractability with respect to the number of free variables.
\begin{theorem}[Partial target solution, standard form] \label{thm:fixed-par-tract}
    Problem \textnormal{(\hyperlink{p:ilfp}{InvLFP})} is fixed-parameter tractable in the parameter $\ell\in\N$ if the LP is given in standard form \eqref{p:lpx-stand} and if $Y = \{\bar y\}\times \R^\ell$ with $\bar y\in\R^{n-\ell}$.
\end{theorem}
\begin{proof}
    We treat the optimistic and pessimistic scenario separately.
    First, we claim that an instance of \eqref{p:ilfpo} with $Y = \{\bar y\}\times \R^\ell$ can be decided by checking if the linear system
    \refstepcounter{equation} \label{eq:fixed-par-tract-opt}
    \begin{alignat}{5}
        & Ay = Bx+b,\quad y\geq 0 \tag{\theequation a} \label{eq:fixed-par-tract-opt-a} \\
        &y_{[n-\ell]} = \bar y \tag{\theequation b} \label{eq:fixed-par-tract-opt-b} \\
        & y_i = 0  && \quad\forall i \in {\mathcal A} \tag{\theequation c} \label{eq:fixed-par-tract-opt-c} \\
        & A^\top\mu \leq Cx+c \tag{\theequation d} \label{eq:fixed-par-tract-opt-d} \\
        & (A^\top\mu)_i = (Cx+c)_i && \quad\forall i \not\in \bar {\mathcal A}\cup {\mathcal A} \tag{\theequation e} \label{eq:fixed-par-tract-opt-e} \\
        & x \in X \tag{\theequation f} \label{eq:fixed-par-tract-opt-f}
    \end{alignat}
    with variables~$x,y,\mu$ is feasible for any subset~${\mathcal A}\subseteq\{n-\ell+1,\dots,n\}$, where we define~$\bar{\mathcal A} \coloneqq  \{i\in[n] \mid \bar y_i=0\}$. Indeed, constraint~\eqref{eq:fixed-par-tract-opt-a} states feasibility of~$y$ and~\eqref{eq:fixed-par-tract-opt-b} is equivalent to membership in~$Y$.
    The optimality of~$y$ is equivalent to the existence of a dual feasible solution~$\mu$ such that~$y$ and~$\mu$ satisfy complementary slackness.
    By the definition of~$\bar{\mathcal A}$, the latter is equivalent to the existence of a partition ${\mathcal A}\subseteq\{n-\ell+1,\dots,n\}$ such that \eqref{eq:fixed-par-tract-opt-c} and \eqref{eq:fixed-par-tract-opt-e} hold.
    Since there exist~$2^\ell$ possible subsets~${\mathcal A}$, we obtain the desired result for the optimistic case.

    The pessimistic scenario can be solved similarly.
    For~$\bar{\mathcal A}$ defined as above and for any subset~$\mathcal A\subseteq\{n-\ell+1,\dots,n\}$, we check whether there exists a parameter~$x\in X$ such that the face~$F$ of~$P(A,Bx+b)$ induced by~$\bar{\mathcal A}\cup\mathcal A$ satisfies the following three conditions:
    \begin{itemize}[align=left,labelwidth={1.25em},labelsep=0.25em]
        \item[(a)] $F$ is the optimal set for~\eqref{p:lpx-stand}.\\[-2.5ex]
        \item[(b)] $F\cap Y\neq\emptyset$.\\[-2.5ex]
        \item[(c)] $F$ is orthogonal to all unit vectors~$\e_h$ for~$h\in[n-\ell]$.
    \end{itemize}
    Assuming~(b), the condition~(c) is equivalent to~$F\subseteq Y$, since~(c) implies that, for all~$h\in[n-\ell]$, all elements of~$F$ agree in the~$h$-th entry and then this entry is~$\bar y_h$ by~(b). Hence, the given instance is a yes-instance if and only if the above test is positive for at least one of the~$2^\ell$ candidate sets~$\mathcal A$.
       
    First note that condition~(c) only depends on~$\mathcal A$, but not on~$x$. Indeed, it is equivalent to the condition that, for all~$h\in[n-\ell]$, the vector~$\e_h$ is spanned by the rows of~$A$ and the unit vectors~$\e_i$, $i\in\bar{\mathcal A}\cup\mathcal A$. This can be verified efficiently.
    So it suffices to check whether both conditions~(a) and~(b) can be satisfied by some common parameter~$x\in X$. Similarly to the previous NP-membership proofs, this is equivalent to the following LP having a positive optimal value:
    \begin{alignat*}{5}
        \max_{x,y,\mu,\varepsilon} & ~~\varepsilon \\
        \st ~~~ & Ay = Bx+b,\quad y\geq 0 \\
        & y_{[n-\ell]} = \bar y  \\
        & y_i = 0 && \qquad\forall i\in {\mathcal A} \\
        & y_i\geq\varepsilon && \qquad\forall i\in\{n-\ell+1,\dots,n\}\setminus {\mathcal A} \\
        & (A^\top\mu)_i = (Cx+c)_i && \qquad\forall i\not\in \bar{\mathcal A}\cup {\mathcal A} \\
        & (A^\top\mu)_i+\varepsilon \leq (Cx+c)_i && \qquad\forall i\in\bar{\mathcal A}\cup {\mathcal A} \\
        & x \in X\;.
    \end{alignat*}
    This concludes the proof.
    
\end{proof}

\begin{remark}[Partial target solution, natural form, OF] \label{rem:xp-natural-of}
    In Theorem~\ref{thm:compl-nat-of} we showed that~(\hyperlink{p:ilfp-S}{InvLFP-S}) is tractable if the LP is given in natural form and only the objective function is parametrized.
    For partial target solutions~{$Y=\{\bar y\}\times\R^\ell$} instead of a target singleton, it is possible to decide~(\hyperlink{p:ilfp}{InvLFP}) for an LP given in natural form and a fixed right-hand side using a similar approach as 
    in Theorem~\ref{thm:fixed-par-tract}.
    In this case, the approach does not yield fixed-parameter tractability, but at least an algorithm that is slicewise polynomial in~$\ell$. 
    The algorithm enumerates every possible subset of active constraints.
    Considering the constraints~$Ay\leq b$ of the LP, the partial target solution leads to a linear system with~$\ell$ free variables.
    Thus, the complete set of active constraints can be characterized by a generating set of at most~$\ell$ active constraints.
    To solve the problem, it then suffices to check for every such possible generating subset~$\mathcal A\subseteq[m]$ with~$|\mathcal A|\leq\ell$ whether it yields a feasible solution. 
    Since there are at most~$\binom{m}{\ell}\leq m^\ell$ such subsets, it follows that~(\hyperlink{p:ilfp}{InvLFP}) with parameter~$\ell$ belongs to~$\text{XP}$ if the LP is given in natural form~\eqref{p:lpx}, if the right-hand side is fixed, and if~$Y=\{\bar y\}\times\R^\ell$ with~$\bar y\in\R^{n-\ell}$.
\end{remark}

Comparing the last results, we have seen that the complexity of~(\hyperlink{p:ilfp}{InvLFP}) with a partially fixed target solution strongly depends on the form of the given~LP:
While the problem is~$\npoly$-complete for any number of fixed variables in natural form (if the right hand side is affected by~$x$), it is fixed-parameter tractable in the number of non-fixed variables for a standard form LP. In particular, in the latter case, it remains tractable if the number of non-fixed variables is bounded. The same is true for natural form LPs if only the objective function is affected by~$x$.

\section{Oracle-defined target set} \label{sec:oracle-target}

In all cases considered so far, the problem~(\hyperlink{p:ilfp}{InvLFP}) belongs to $\npoly$. In the optimistic scenario, it is easy to see that the complexity does not even exceed~$\npoly$ if instead of a polyhedron~$Y$ we consider only its integer points~$Y\cap\Z^n$ as target set.
Indeed, we can simply add integrality constraints~$y\in\Z^n$ to the linear system~\eqref{eq:polyh-np-opt}.  
Since we only need a feasibility certificate for the resulting mixed-integer linear system, this restriction does not increase the complexity.

Simply put, and more generally, for the verification of optimistic yes-instances it suffices to verify independently that a solution~$y$ is optimal for an LP with some parameter~$x\in X$ and that~$y$ belongs to~$Y$. 
This suggests that it might be possible to generalize Theorem~\ref{thm:polyh-np} even further and require only some kind of membership oracle for~$Y$.
However, it is easy to construct a yes-instance of~\eqref{p:ilfpo}, with a non-polyhedral target set~$Y$, such that the unique parameter~$x$ making an element of~$Y$ optimal is irrational, so that~$x$ itself cannot be chosen as certificate. The solution is to consider a weaker type of oracle.

\begin{definition}\label{def:out-mem-oracle}
    A \emph{weak outer membership oracle} for~$Y\subseteq\R^n$ receives as input a point $y\in\R^n$ and some~$\delta>0$ and either asserts $y\not\in S(Y,\delta)$ or~$y\in S(Y,2\delta)$.
\end{definition}
\noindent Here we define $S(Y,\varepsilon)\coloneqq \{x\in\R^n\mid \exists y\in Y\colon ||x-y||\le\varepsilon\}$. Definition~\ref{def:out-mem-oracle} is closely related to the common definition of a weak membership oracle, as given in Grötschel et al.~\cite{grotschel1993geometricalgorithms}, but applied to a neighborhood of~$Y$ instead of~$Y$ itself.

\begin{theorem}[Oracle-based target set, optimistic scenario] \label{thm:weak-np-memb-oracle}
    Let $Y\subseteq\R^n$ be bounded and given by a weak outer membership oracle. Then, if an instance of~\eqref{p:ilfpo} with target set~$Y$ is a yes-instance, it is possible to verify for any~$\delta > 0$ that the same instance with target set~$S(Y,\delta)$ is a yes-instance, in time polynomial in the input size, in $1/\delta$, and in the bounding radius of~$Y$.
\end{theorem}
\begin{proof}
    Let a yes-instance of \eqref{p:ilfpo} with target set~$Y$ be given and assume that~$\delta>0$.
    Then the same instance with target set $S(Y,\delta)$ is a yes-instance as well.
    We will prove that there exists a solution~$x\in X$ to the latter instance that has polynomial encoding length.
    By assumption, there exist~$\bar x\in X$ and~$\bar y\in Y$ such that~$\bar y\in \bar F$, where
    $$\bar F\coloneqq \arg\min\{ (C\bar x+c)^\top y \mid Ay \leq B\bar x + b\}\;.$$
    The first part of the desired certificate for the modified instance will be the set of active indices
    $$\mathcal A \coloneqq  \{j\in[m] \mid (Ay)_j = (B\bar x+b)_j\ \forall y\in\bar F\}$$
    for~$\bar F$ and~$\bar x$.
    Next, since~$\bar y$ is bounded, there exists a point~$\hat y\in B(\bar y,n^{-1/2}\cdot\delta/4)$ of polynomial encoding length in the input size, in $1/\delta$, and in the bounding radius~$R$ of~$Y$, which will be the second part of the certificate. In particular, we have~$\hat y\in B(\bar y,\delta/4)$, so that the weak outer membership oracle applied to~$\hat y$ and~$\delta /4$ can only return the answer~$\hat y\in S(Y,\delta/2)$.
    In summary, the modified instance can be verified as a yes-instance by calling the oracle for $\hat y$ and~$\delta/4$ and by showing that the following linear system in the variables $x,y,\mu$ is feasible:
    \refstepcounter{equation} \label{eq:weak-np-memb-oracle-verif}
    \begin{alignat}{5}
        & Ay \leq Bx+b \tag{\theequation a} \label{eq:weak-np-memb-oracle-verif-a} \\
        & (Ay)_j = (Bx+b)_j \qquad\qquad                           && \forall j\in \mathcal A \tag{\theequation b} \label{eq:weak-np-memb-oracle-verif-b} \\
        & A^\top \mu = Cx+c,\ \mu \leq 0 \tag{\theequation c} \label{eq:weak-np-memb-oracle-verif-c} \\
        & \mu_j = 0\quad                                  && \forall j\not\in \mathcal A \tag{\theequation d} \label{eq:weak-np-memb-oracle-verif-d} \\
        &y_i - \hat y_i \leq n^{-\frac 12} \delta/2\qquad   && \forall i\in[n] \tag{\theequation e} \label{eq:weak-np-memb-oracle-verif-e1} \\
        &\hat y_i-y_i\leq n^{-\frac 12} \delta/2\qquad   && \forall i\in[n] \tag{\theequation f} \label{eq:weak-np-memb-oracle-verif-e2} \\        & x\in X. \tag{\theequation g} \label{eq:weak-np-memb-oracle-verif-f}
    \end{alignat}
    Indeed, by assumption, there exists a $\bar\mu\in\R^m$ such that $(\bar x, \bar y,\bar\mu)$ solves the system~\eqref{eq:weak-np-memb-oracle-verif}, where we use $||y-\hat y||_\infty \le ||y-\hat y||_2\le n^{-1/2}\delta/4$.
    
    Conversely, given the certificate~$\mathcal A$ and~$\hat y$, we claim that if our oracle call for~$\hat y$ and~$\delta/4$ returns the answer~$\hat y\in S(Y,\delta/2)$ and if~$(x,y,\mu)$ solves the above linear system, then $(x,y)$ is a solution to the modified instance. In fact, the constraints~\eqref{eq:weak-np-memb-oracle-verif-a}--\eqref{eq:weak-np-memb-oracle-verif-d} ensure that $y$ is feasible and optimal for the LP with parameter $x$.
    From \eqref{eq:weak-np-memb-oracle-verif-e1} and \eqref{eq:weak-np-memb-oracle-verif-e2} it follows that~\mbox{$\|y-\hat y\|_2 \leq n^{1/2}||y-\hat y||_\infty\le\delta/2$}.
    Since~$\hat y\in S(Y,\delta/2)$, this implies $y\in S(Y,\delta)$, verifying that $(x,y)$ is indeed a solution to the modified instance.
    
\end{proof}

We next consider the pessimistic case. Here, the main additional challenge, compared to the optimistic case, is to ensure that the set of optimal solutions to the LP parametrized by~$x$ is entirely contained in the target set~$Y$.
Especially for the case where~$Y$ is a general polyhedron, this approach was harder to realize than in the optimistic version, and the idea relied heavily on the polyhedral structure of~$Y$; see the proof of Theorem~\ref{thm:polyh-np}.
This indicates that a generalization with regard to the target set in the pessimistic version is not as straightforward as in the optimistic version. 
In fact, the following results prove that a similar approach as for the optimistic version is not possible at all.
\begin{lemma}\label{lem:oracle}
    If~$Y\subseteq\R^n$ is given by a strong membership oracle, then it is not possible to verify, using a polynomial number of oracle calls, that $[0,1]^n\subseteq Y$, even if~$Y$ is a polynomially encodable polytope.
\end{lemma}
\begin{proof}
    Let~$V=\{v_1,\dots,v_{2^n}\}$ be the set of all vertices of~$[0,1]^n$ and consider the polytopes~$Y_1,\dots,Y_{2^n}$ defined by $Y_i \coloneqq  \conv\left(V\setminus \{v_i\}\right)$.
    To verify $[0,1]^n\subseteq Y$, we must rule out that $Y=Y_i$ for any $i=1,\dots,2^n$ during the verification of the certificate.
    However, we can only gain information about~$Y$ by calling the strong membership oracle, and 
    each call of the oracle can rule out at most one of the polytopes $Y_1,\dots,Y_{2^n}$.
    Hence, after polynomially many oracle calls, for at least one~$i\in[2^n]$,
    the verification must yield the same result for~$[0,1]^n$ and for~$Y_i$, showing that a verification is not yet possible.
    Since~$[0,1]^n$ as well as every polytope~$Y_i$ can be defined by at most~$2n+1$ linear inequalities, this also proves the last claim.
    
\end{proof}

\begin{theorem}[Oracle-based target set, pessimistic scenario] \label{thm:pess-no-verif}
    If the target set~$Y$ is given by a strong membership oracle, then it is not possible to define certificates for yes-instances of~\eqref{p:ilfpp} that can be verified in polynomial time.
    This even holds if~$Y$ is a polynomially encodable polytope.
\end{theorem}
\begin{proof}
    This follows from Lemma~\ref{lem:oracle} by considering the instance of~\eqref{p:ilfpp} that consists of the parameter set~$X=\{0\}$ and the LP
    \begin{equation*}
        \begin{aligned}
            \min\quad&      0^\top y \\
            \st\quad&       y\in[0,1]^n\;,
        \end{aligned}
    \end{equation*}
    which is actually independent of the parameter~$x$.
    Indeed, the optimal set of this~LP is~$[0,1]^n$, so that we have to verify~$[0,1]^n\subseteq Y$ as in Lemma~\ref{lem:oracle}.
    
\end{proof}
The statements of Theorem~\ref{thm:weak-np-memb-oracle} and Theorem~\ref{thm:pess-no-verif} assume different types of oracles. However, it is obvious that the negative result of the latter holds in particular when, instead of a strong membership oracle, only a weak outer membership oracle is given. In fact, the negative results of Lemma~\ref{lem:oracle} and Theorem~\ref{thm:pess-no-verif} even remain true when~$Y$ is given by a linear optimization oracle. This follows from the fact that a linear optimization oracle cannot distinguish more efficiently between~$[0,1]^n$ and the sets~$Y_i$ than a membership oracle.

\section{Conclusion} \label{sec:conclusion}

In this paper, we have presented a comprehensive analysis of the inverse linear feasibility problem~(\hyperlink{p:ilfp}{InvLFP}). We showed that its complexity depends on several aspects of the problem, in particular on the type of target set, the form of the~LP, and the optimistic or pessimistic scenario.
For $\npoly$-hardness results, we used reductions from 3-SAT, building on established $\npoly$-hardness proofs of bilevel linear programming.
The tractability results rely on linearization of (unique) optimality conditions based on (strict) complementary slackness.
The main complexity results are summarized in Table~\ref{tab:overview}.

\begin{table}[h]
    \centering
    \begin{adjustbox}{max width=\textwidth}
    \begin{threeparttable}
        \begin{tabular}{llllll} 
        {\bf Target set} & {\bf LP form}     & {\bf Param.}    & {\bf Scen.} & {\bf Result} & {\bf Reference} \\ \toprule
        Singleton   & natural  & RHS   &       & $\npoly$-hard & Th.~\ref{thm:compl-nat} \\
                    & natural  & OF    &       & tractable & Th.~\ref{thm:compl-nat-of} \\
                    & standard &       &       & tractable & Th.~\ref{thm:comp-stand} \\ \midrule
        Basis       & standard & RHS   & opt.  & $\npoly$-hard & Th.~\ref{thm:basis-stand-opt} \\
                    & standard & OF    & opt.  & tractable  & Rem.~\ref{rem:basis-stand-opt-of} \\
                    & standard &       & pess. & tractable & Th.~\ref{thm:basis-stand-pess} \\\midrule
        Polyhedron \\
        \quad $\{y_i=0\}$ &         &       &       & $\npoly$-hard & Th.~\ref{thm:one-variable-nat-rhs}, Th.~\ref{thm:one-variable-nat-of} \\
        \quad $\{\bar y\}\times\R^\ell$  & standard    &     &  & FPT in $\ell$ & Th.~\ref{thm:fixed-par-tract} \\ 
        \quad $\{\bar y\}\times\R^\ell$  & natural     & RHS  &  & $\npoly$-hard for fixed $\ell$ & Th.~\ref{thm:compl-nat} \\ 
        \quad $\{\bar y\}\times\R^\ell$  & natural     & OF  &  & $\text{XP}$ in $\ell$ & Rem.~\ref{rem:xp-natural-of} \\ \midrule 
        Oracle based &            &       & opt.  & ``weakly verifiable'' & Th.~\ref{thm:weak-np-memb-oracle} \\
                        &           &         & pess. & not verifiable & Th.~\ref{thm:pess-no-verif} \\ \bottomrule
    \end{tabular}
        \caption{Complexity overview of (InvLFP). The table presents the main 
        complexity results for (InvLFP) shown in this paper, depending on the structure of the target set, the form of the underlying LP, which parameters in the LP can be modified, and the type of scenario. Empty cells in the columns for the form, parameter, and scenario indicate that the result holds for any option for this characteristic.}
    \label{tab:overview}
    \end{threeparttable}
    \end{adjustbox}
\end{table}

For the optimistic scenario, we showed that~\eqref{p:ilfpo} with a single target solution is $\npoly$-hard if the LP is given in natural form and the right-hand side is parametrized.
This changes if the LP is given in standard form, since a standard form target solution provides enough information to linearize the optimality conditions, rendering it tractable.
However, if the target solution is replaced by a target basis, the resulting problem becomes $\npoly$-complete again, because a potential degeneracy of the basis leads to exponentially many possibilities to satisfy the optimality conditions.
For partial target solutions, we proved that~\eqref{p:ilfpo} with a standard form LP admits an FPT-algorithm in the number of unfixed variables of the target solution, but becomes $\npoly$-hard if that number is unbounded and at least one target variable is fixed.
For a natural form LP with a partial target solution, the problem admits a slicewise polynomial time algorithm in the number of unfixed variables in the target solution.

For the pessimistic scenario, we showed the same complexity results in most cases, with slightly different methods and a few notable exceptions.
In the case of a single target solution, the complexity is equal to the optimistic scenario, but when switching to a target basis, the resulting problem remains tractable.
Roughly speaking, the reason is that unique optimality requires the absence of degeneracy, which was the main source of hardness in the optimistic case.
For polyhedral target sets, we showed the same hardness results for the pessimistic and the optimistic case.
In both scenarios, we obtained an upper bound for the complexity of~(\hyperlink{p:ilfp}{InvLFP}) by showing membership in $\npoly$.
However, the two proofs contained significant differences, which become apparent when trying to extend this upper bound to general target sets based on membership oracles:
For~\eqref{p:ilfpo}, the~$\npoly$-membership only requires to solve an LP and check membership in the target set, up to some rounding that is necessary in the case of non-polyhedral target sets.
For~\eqref{p:ilfpp}, we proved that a similar approach is not possible, as even deciding whether a simple target polytope is a subset of another polytope is not possible in polynomial time if the latter is only given via a membership oracle.
This emphasizes the additional challenge of the pessimistic case with respect to the optimistic case, namely that it includes a condition concerning the inclusion of a polyhedron (the optimal set) in~$Y$. Also in bilevel optimization, this often creates additional difficulties; see, e.g.,~\cite{Buchheim2023BilevelInNP}.

It would be interesting to determine how the complexity of the considered problems behaves if the parametric LP is replaced by a harder parametrized problem.
Note that even though~$\npoly$-membership is preserved for mixed-integer target sets, as mentioned in the beginning of Section~\ref{sec:oracle-target}, it is easy to show that considering a parametric integer linear program instead of an LP already leads to a~$\Sigma^\poly_2$-hard problem in general.
For a related problem with a parametric ILP where the parameter and target set are also integer, Kurtz et al.~\cite{Kurtz2025CounterfactualILP} show~$\Sigma_2^\poly$-completeness if the constraint matrix~$A$ is affinely parametrized.
However, it is not clear whether all our hardness results are simply lifted onto the next level of the polynomial hierarchy in this case. This is left as future work.

\pagebreak

\bibliographystyle{abbrvnat}
\bibliography{pilfp-bibliography}

\end{document}